\documentclass[12pt]{amsart}
\usepackage{amsmath,amssymb,graphicx}
\usepackage[all]{xy}
\usepackage[bbgreekl]{mathbbol}
\newtheorem{thm}{Theorem}[section]
\newtheorem{cor}[thm]{Corollary}

\newtheorem{lem}[thm]{Lemma}
\newtheorem{prop}[thm]{Proposition}
\theoremstyle{definition}
\newtheorem{defn}[thm]{Definition}
\theoremstyle{remark}
\newtheorem{rem}[thm]{Remark}
\numberwithin{equation}{section}
\DeclareFontFamily{U}{rsf}{} \DeclareFontShape{U}{rsf}{m}{n}{
  <5> <6> rsfs5 <7> <8> <9> rsfs7 <10->  rsfs10}{}
\DeclareMathAlphabet{\mathscr}{U}{rsf}{m}{n}
\newcommand{\mycal}[1]{\mathscr{#1}}

\newcommand{\mto}[1]{\stackrel{#1}{\longrightarrow}}
\newcommand{\abs}[1]{\left\vert#1\right\vert}

\newcommand{\A}{\mathcal{A}}
\newcommand{\B}{\mathcal{B}}

\newcommand{\Pc}{{\mathcal P}}
\newcommand{\Qc}{{\mathcal Q}}
\newcommand{\Zb}{{\mathbb Z}}

\newcommand{\Xc}{\mathcal X}
\newcommand{\Yc}{\mathcal Y}
\renewcommand{\imath}{\sqrt{-1}}
\newcommand{\Sc}{\mycal S}

\newcommand{\kk}{\mathbb{k}}

\newcommand{\mem}{\boldsymbol{ ^{{\small\backprime}}\!\!\!\cap\!\!\!\!
\_\,}}
\newcommand{\bbb}{\Box\!\!\!\!\Box}
\newcommand{\Ab}{\mathbb A}
\newcommand{\Eb}{\mathbb E}

\newcommand{\cf}{2\pi\imath\,}

\newcommand{\np}{ }
\newcommand{\dbar}{\overline{\partial}}
 \DeclareMathOperator{\Hom}{Hom}
\DeclareMathOperator{\Ext}{Ext} \DeclareMathOperator{\End}{End}

\makeindex
\begin{document}

\title[$\mem$-Branes]{ Duality and equivalence of module categories in noncommutative
geometry II: \\
Mukai duality for holomorphic noncommutative tori}
\author{Jonathan Block}%
\address{Department of Mathematics, University of Pennsylvania, Philadelphia, PA 19104}%
\email{blockj@math.upenn.edu}%

\thanks{J.B. partially supported by NSF grant DMS02-04558}%
\dedicatory{To my parents}

\begin{abstract}
This is the second in a series of papers intended to set up a
framework to study categories of modules in the context of
non-commutative geometries. In \cite{mem} we introduced the basic
DG category $\Pc_{\A^\bullet}$, the perfect category of
$\A^\bullet$,  which corresponded to the category of coherent
sheaves on a complex manifold. In this paper we enlarge this
category to include objects which correspond to quasi-coherent
sheaves. We then apply this framework to proving an equivalence of
categories between  derived categories on the noncommutative
complex torus and on a holomorphic gerbe on the dual complex
torus.

\end{abstract}

 \maketitle \tableofcontents
\newpage\section{Introduction}
This is the second in a series of papers. Its purpose is twofold.
First, we continue to set up a framework to study categories of
modules in the context of noncommutative geometries, which are
integrable in some sense. This integrability we encode in a
(curved) differential graded algebra $(\A^\bullet, d, c)$. The
second is to apply our framework to an interesting non-trivial
example. In particular, we prove an extension of Mukai duality
from the classical case of complex tori, to that of noncommutative
tori.

 In \cite{mem} we introduced the basic DG category
$\Pc_{\A^\bullet}$, the \index{perfect category} of $\A^\bullet$,
which corresponded to the category of coherent sheaves on a
complex manifold.

One of the purposes of this paper is to develop more of the
apparatus. In particular, we show how to get some of the six
operations of Grothendieck in our context. They are all some
version of tensoring  with a twisted bimodule. However, in order
to do this, we need to enlarge our basic perfect category of
modules $\Pc_{\A^\bullet}$ to a larger category
$q\Pc_{\A_\bullet}$, the quasi-perfect category of $\A^\bullet$.
This corresponds to the category of quasi-coherent sheaves on a
complex manifold. The category of quasi-coherent sheaves on a
complex manifold is less well known for complex manifolds than it
is for algebraic schemes.

There are at least two reasons for enlarging the category. Most
pressingly, the derived pushforward  is most naturally defined
into the quasi-perfect category. For example, in complex geometry,
to define the derived pushforward, one performs the following
steps:

\begin{enumerate}
\item Resolve by injective sheaves, which can only be done in the
quasi-coherent category, \item Push forward. This remains a
quasi-coherent complex of sheaves. \item Prove that the direct
image is equivalent to a complex of coherent sheaves, i.e.
Grauert's direct image theorem.
\end{enumerate}

If we don't care to end up back in the coherent category, we can
stop after step (2). This is what we do in this paper. Finiteness
conditions and Grauert's direct image theorem are the subject of
the third paper in the series, \cite{samech}.

The second reason for introducing this larger category is that on
a complex manifold, there can potentially be very few coherent
sheaves. Certainly, not nearly enough to determine the manifold up
to isomorphism. Quasi-coherent sheaves provides a larger category
where more invariants of complex manifolds might be found.

Most of the paper is concerned with an application of our
framework: to formulate and prove a deformed version of Mukai
duality,  which we now explain. let $X$ be a complex torus. Thus
$X=V/\Lambda$ where $V$ is a $g$-dimensional complex vector space
and $\Lambda\cong \Zb^{2g}$ is a lattice in $V$. Let $X^\vee$
denote the dual complex torus. This can be described in a number
of ways:

\begin{enumerate}
\item as $\mbox{Pic}^0(X)$, the variety of holomorphic line
bundles on $X$ with first Chern class $0$ (i.e., they are
topologically trivial); \item as the moduli space of flat unitary
line bundles on $X$. This  is the same as the space of unitary
representations of $\pi_1(X)$, but it has a complex structure that
depends on that of $X$; \item and most explicitly as
$\overline{V}^\vee/\Lambda^\vee$ where $\Lambda^\vee$ is the dual
lattice,
$$\Lambda^\vee=\{v\in\overline{V}^\vee \, |\, \mbox{ Im } <v,\lambda>\in
\Zb\,\, \forall \lambda \in\Lambda\}.$$ Here $\overline{V}^\vee$
consists of conjugate linear homomorphisms from $V$ to
$\mathbb{C}$.
\end{enumerate}

We now describe Mukai duality. On $X\times X^\vee$ there is a
canonical line bundle, $\Pc$, the Poincar\'{e} bundle which is
uniquely determined by the following universal  properties:
\begin{enumerate}
\item $\Pc|X\times\{p\}\cong p$ where $p\in X^\vee$ and is
therefore a line bundle on $X$. \item $\Pc|\{0\}\times X^\vee$ is
trivial.
\end{enumerate}
Now Mukai duality says that there is an equivalence of derived
categories of coherent sheaves
$$D^b(X)\to D^b(X^\vee)$$
induced by the functor
$${\mathcal F}\mapsto p_{2*}(p_1^*{\mathcal F}\otimes \Pc)$$
where $p_i$ are the two obvious projections. The relation of this
statement with the Baum-Connes conjecture is discussed in
\cite{mem}.

Besides the usual deformations of a complex manifold $X$ as a
complexmanifold, there are interesting extra deformations that are
derived from the philosophy expounded by Bondal, Drinfeld and
Kontsevich: The most general way to deform a space $X$ is by
deforming the derived category $D^{b}(X)$ of sheaves or a DG
enhancement (such as $\Pc_{\A^{0,\bullet}(X)}$) as a DG category.
The infinitesimal deformations of $D^{b}(X)$ are given by the
second Hochschild cohomology of $X$:
\[
HH^{2}(X) := \Ext^{2}_{X\times X}({\mathcal
  O}_{\Delta},{\mathcal O}_{\Delta}).
\]
There is a `Hodge type' decomposition:
\[
HH^{2}(X) = H^{0}(X;\wedge^{2}T_{X})\oplus H^{1}(X;T_{X}) \oplus
H^{2}(X;{\mathcal O}_{X})
\]
\begin{enumerate}
\item The term  $H^{1}(X;T_{X})$ corresponds to the classical
deformations of $X$ as a complex manifold. \item
$H^{0}(X;\wedge^{2}T_{X})$ consists of global holomorphic Poisson
structures and correspond to deformations of $X$ as a
noncommutative space. \item The most mysterious term
$H^{2}(X;{\mathcal O}_{X})$ corresponds to deformations of the
trivial ${\mathcal O}^\times$ gerbe to a non-trivial ${\mathcal
O}^\times$-gerbe.
\end{enumerate}

In the case above, where $X=V/\Lambda$ is a complex torus and
$X^\vee=\overline{V}^\vee/\Lambda^\vee$ its dual torus, the
equivalence of categories implemented by the Poincar\'{e} bundle
establishes the following identification of the terms of the
Hochschild cohomology:
\begin{equation}\label{eq-def-corr}
\begin{array}{lllll}
H^0(X;\wedge^{2}T_{X}) & \cong & \wedge^2 V &\cong &
H^2(X^\vee;\mathcal{O}) \\
H^1(X;T_X) & \cong & V\otimes \overline{V}^\vee &\cong &
H^1(X^\vee;T_{X^\vee})\\
H^2(X;\mathcal{O}) & \cong & \wedge^2 \overline{V}^\vee &\cong &
H^2(X^\vee;\wedge^2 T_{X^\vee})
\end{array}\end{equation}

\noindent What this suggests, is that if we deform $X$ to a
holomorphic noncommutative torus, then the dual will deform to a
holomorphic gerby torus, and vice versa. This is where one should
look for a derived equivalence of categories.

This derived equivalence  was carried out in the context of
deformation quantizations by Ben-Bassat, Block and Pantev,
\cite{BBP}. Here we will prove this derived equivalence in our
context which has the advantage over \cite{BBP} of applying to
informal (i.e. non-formal) deformations.

The author would like to thank Oren Ben-Bassat, Calder Daenzer,
Anton Kapustin, Tony Pantev and Betrand Toen for conversations and
suggestions during this work.

\section{The quasi-perfect category of modules over a differential graded algebra}
In this section we make our extension of the perfect category of a
(curved) DGA.

\subsection{Review of the perfect DG category of a
curved DGA}

Here we briefly review the DG-category that we assigned to a
curved DGA. See \cite{mem} for more details. Also, see \cite{BK},
\cite{Ke}, \cite{Dr} for facts about DG-categories.
\begin{defn} A curved DGA \cite{PP} (Schwarz calls them $Q$-algebras) is a triple
$(\A^\bullet, d, c)$ where $\A^\bullet$ is a (non-negatively)
graded algebra over $k$ with a derivation
\[d:\A^\bullet\to \A^{\bullet+1}\]
which satisfies the usual Leibnitz relation but
\[d^2(a)=[c,a]\]
where $c\in \A^2$ is a fixed element (the curvature). Furthermore
we require the Bianchi identity $dc=0$. \end{defn} \noindent We
will always set $\A=\A^0$.

A DGA is the special case where $c=0$. Note that $c$ is part of
the data and even if $d^2=0$, that $c$ might not be $0$, and gives
a non DGA example of a curved DGA. This is in fact the case for
our gerby complex torus described in the introduction. The
prototypical example of a curved DGA is $(\A^\bullet(M,\End(E)),
Ad \nabla, F)$ of differential forms on a manifold with values in
the endomorphisms of a vector bundle $E$ with connection
${\nabla}$ and curvature $F$.

\np Let $(\A^\bullet,d, c)$ be a curved DGA. Let $E^\bullet$ be a
$\Zb$-graded (bounded in both directions) right $\A$-module which
is finitely generated and projective.

\begin{defn}\label{Z-connection} A
$\Zb$-connection $\Ab$  is a $k$-linear  map
\[\Ab: E^\bullet\otimes_\A\A^\bullet\to E^\bullet\otimes_\A\A^\bullet
\]
of total degree one, which satisfies the usual Leibnitz condition
\[  \Ab(e\omega)=(\Ab(e))\omega+(-1)^e e d\omega\]
\end{defn}
Such a connection is determined by its value on $E^\bullet$. Let
$\Ab^k$ be the component of $\Ab$ such that $\Ab^k:E^\bullet\to
E^{\bullet-k+1}\otimes_\A\A^k$, thus
$\Ab=\Ab^0+\Ab^1+\Ab^2+\cdots$.  It is clear that $\Ab^1$ is a
connection on $E^\bullet$ in the ordinary sense and that $\Ab^k$
is $\A$-linear for $k\ne 1$.

\begin{defn} For a curved DGA $(\A^\bullet,d,c)$, an object in the {\em
perfect DG-category} $\Pc_{\A^\bullet}$ (called a perfect twisted
complex, or a module) is a $\Zb$-graded (but bounded in both
directions) right module $E^\bullet$ over $\A$ which is finitely
generated and projective, a $\Zb$-connection
\[\Eb: E^\bullet\otimes_\A\A^\bullet\to E^\bullet\otimes_\A\A^\bullet
\]
that satisfies the integrability condition
\[\Eb\circ\Eb(e)=-e\cdot c\]
The minus appears because we are dealing with right modules.

\np The morphisms between two such objects $E_1=(E_1^\bullet,
\Eb_1)$ and $E_2=(E_2^\bullet,\Eb_2)$ of degree $k$ are
\[
\Hom^k_{\Pc_\A^\bullet}(E_1,E_2)=\{\phi:E_1^\bullet\otimes_\A
\A^\bullet\to E_2^{\bullet}\otimes_\A \A^\bullet\,\,|\,\,
\phi(ea)=(-1)^{k|a|}\phi(e)a\}
\]
with differential defined in the standard way
\[
d(\phi)(e)= \Eb_2(\phi(e))-(-1)^{\abs{\phi}}\phi(\Eb_1(e))
\]
\end{defn}

Again, such a $\phi$ is determined by its restriction to
$E_1^\bullet$ and if necessary we denote the component of $\phi$
that maps \begin{equation}\label{pcomp} E_1^\bullet\to
E_2^{\bullet+k-j}\otimes_\A \A^j \end{equation} by $\phi^j$.

\begin{prop} For a curved DGA $(\A^\bullet, d, c)$ the category $\Pc_{\A^\bullet}$ is a DG-category.
\end{prop}
\np This is clear from the following lemma.

\begin{lem}
Let $E_1$, $E_2$ be modules over the curved DGA
$(\A^\bullet,d,c)$. Then the differential defined above
\[d:\Hom^\bullet_{\Pc_{\A^\bullet}}(E_1,E_2)\to \Hom^{\bullet+1}_{\Pc_{\A^\bullet}}(E_1,E_2)\]
satisfies $d^2=0$.
\end{lem}
\proof This is a simple check, and follows because the curvature
terms from $E_1$ and $E_2$ cancel.
\endproof

\begin{defn} A morphism $f:X\to Y$ in $\Pc_{\A^\bullet}$ which is
closed and of degree $0$ is a quasi-isomorphism if and only if for
all objects $A$ of $\Pc_{\A^\bullet}$ the morphism induced by
post-composing with $f$,
\[
f^\#:\Hom_{\Pc_{\A^\bullet}}^\bullet(A,X)\to
\Hom_{\Pc_{\A^\bullet}}^\bullet (A,Y)
\]
is a quasi-isomorphism of complexes. \end{defn} In the case of the
Dolbeault algebra, it is classical that this recovers the usual
notion of a morphism of complexes of holomorphic vector bundles
being a quasi-isomorphism.

In \cite{mem} we proved the following criterion for being a
quasi-isomorphism.
\begin{prop}
Under the additional assumption that each $\A^p$ is flat as an
$\A$-module, a closed morphism $\phi\in
\Hom^0_{\Pc_{\A^\bullet}}(E_1,E_2)$ is a quasi-isomorphism if
$\phi^0:(E_1^\bullet,\Eb_1^0)\to (E_2^\bullet,\Eb_2^0)$ is a
quasi-isomorphism of complexes of $\A$-modules.
\end{prop}

We also discussed how to construct functors between categories of
the form $\Pc_{\A^\bullet}$. Let $(\A_1^\bullet,d,c_1)$ and
$(\A^\bullet_2,d,c_2)$ be two curved DGAs. Consider the following
data, $\mathcal{X}=(X^\bullet,\mathbb{X})$ where
\begin{enumerate}
\item $X^\bullet$ is a graded finitely generated projective
right-$\A_2$-module, \item $\mathbb{X}:X^\bullet\to
X^\bullet\otimes_{\A_2}\A_2^\bullet$ is a $\Zb$-connection, \item
$\A^\bullet_1$ acts on the left of
$X^\bullet\otimes_{\A_2}\A_2^\bullet$ satisfying

\[
a\cdot (x\cdot b)=(a\cdot x)\cdot b
\]
and
\[
\mathbb{X}(a\cdot(x\otimes b))=da\cdot(x\otimes
b)+a\cdot\mathbb{X}(x\otimes b)
\]
for $a\in \A_1^\bullet$, $x\in X^\bullet$ and $b\in \A_2^\bullet$,
 \item
$\mathbb{X}$ satisfies the following condition
\[
\mathbb{X}\circ\mathbb{X}(x\otimes b)=c_1\cdot (x\otimes
b)-(x\otimes b)\cdot c_2
\]
 on the complex
$X^\bullet\otimes_{\A_2}\A_2^\bullet$.

\end{enumerate}
Let us call such a pair $\mathcal{X}=(X^\bullet,\mathbb{X})$ an
$\A_1^\bullet-\A_2^\bullet$-twisted  bimodule.

\np Given a $\A_1^\bullet-\A_2^\bullet$-twisted  bimodule
$\mathcal{X}=(X^\bullet,\mathbb{X})$, we can then define a
DG-functor
\[ \mathcal{X}_*:\Pc_{\A_1^\bullet}\to \Pc_{\A_2^\bullet}\]
by $\mathcal{X}_*(E^\bullet, \Eb)$ is the twisted complex
\[(E^\bullet\otimes_{\A_1}X^\bullet,\Eb_2)\]
where $\Eb_2(e\otimes x)=\Eb(e)\cdot x+e\otimes \mathbb{X}(x)$,
where the $\cdot$ denotes the action of $\A_1^\bullet$ on
$X^\bullet\otimes_{\A_1}\A_2^\bullet$.  One easily checks that
$\mathcal{X}_*(E)$ is an object of $\Pc_{\A_2^\bullet}$. We will
write $\mathbb{E}\# \mathbb{X}$ for $\mathbb{E}_2$.

\begin{enumerate}
\item Given an $\A_1^\bullet-\A_2^\bullet$-twisted  bimodule
$\mathcal{X}=(X^\bullet,\mathbb{X})$ and an
$\A_2^\bullet-\A_3^\bullet$-twisted  bimodule
$\mathcal{Y}=(Y^\bullet,\mathbb{Y})$ we can form an
$\A_1^\bullet-\A_3^\bullet$-twisted  bimodule
\[
\mathcal{X}\otimes_{\A_2}\mathcal{Y}=(X^\bullet\otimes_{\A_2}Y^\bullet,\mathbb{X}\#\mathbb{Y})
\]
in the same way that we defined the functor $\mathcal{X}_*$.
Moreover, it is clear that the functors $\mathcal{Y}_*\circ
\mathcal{X}_*$ and $(\mathcal{X}\otimes_{\A_2}\mathcal{Y})_*$ are
naturally isomorphic. \item Given two curved DGA's,
$(\A_1^\bullet,d_1,c_1)$ and $(\A_2^\bullet,d_2,c_2)$ a
homomorphism from $\A_1^\bullet$ to $\A_2^\bullet$ is a pair
$(f,\omega)$ where $f:\A_1^\bullet\to\A_2^\bullet$ is a morphism
of graded algebras, $\omega\in \A_2^1$ and they satisfy
\begin{enumerate}
\item $f(d_1a_1)=d_2f(a_1)+[\omega,f(a_1)]$ and \item
$f(c_1)=c_2+d_2\omega+\omega^2$.
\end{enumerate} To such a homomorphism we associate the twisted
bimodule $\mathcal{X}_f$ where  $X_f^\bullet=\A_2$ in degree $0$.
$\A_1^\bullet $ acts by the morphism $f$ and the $\Zb$-connection
is
\[ \mathbb{X}_f(a_2)=d_2(a_2)+\omega\cdot a_2.\]
\item As a special case of the previous example, when $\phi:X\to
Y$ is a holomorphic map between complex manifolds and
$f=\phi^*:\A^{0,\bullet}(Y)\to \A^{0.\bullet}(X)$ is the induced
map on the Dolbeault DGAs, then the DG-functor
\[\mathcal{X}_{f*}:\Pc_{\A^{0,\bullet}(Y)}\to
\Pc_{\A^{0.\bullet}(X)}
\]
is just the pullback functor of coherent sheaves.
\end{enumerate}

\subsection{The Quasi-perfect category}
The need to define analogues of derived pushforwards of coherent
sheaves drives us to introduce the larger quasi-perfect category
$q\Pc$. Let $X$ and $Y$ be complex manifolds. Let $p_1$ (resp.
$p_2$) denote the projection from $X\times Y$ to $X$ (resp. $Y$).
We will start by describing the twisted bimodule that should
implement the pushforward
\[
p_{2*}:\Pc_{\A^{0,\bullet}(X\times Y)}\to \Pc_{\A^{0.\bullet}(Y)}
\]
Let \[M^\bullet=\Gamma(X\times Y;\wedge^\bullet p_1^* T^{0,1}X)\]
be the space of smooth forms along the fiber of the projection
$p_2$. This has a natural right action of $\A(Y)$. Define a
$\mathbb{Z}$-connection
\begin{equation}\label{Z-conn}
\mathbb{M}:M^\bullet\to
M^\bullet\otimes_{\A(Y)}\A^{0,\bullet}(Y)\cong
\A^{0,\bullet}(X\times Y)
\end{equation}
by $\mathbb{M}=\mathbb{M}^0+\mathbb{M}^1$ where
$\mathbb{M}^0=\dbar_X$ is the $\dbar$-operator along the fiber and
$\mathbb{M}^1=\dbar_Y$. Note that $\mathbb{M}^0$ is
$\A(Y)$-linear. To complete the construction of a twisted bimodule
we need to construct an action of $\A^{0,\bullet}(X\times Y)$ on
$M^\bullet\otimes_{\A(Y)}\A^{0,\bullet}(Y)$. Using the isomorphism
in \eqref{Z-conn} this is just left multiplication of
$\A^{0,\bullet}(X\times Y)$ on itself. All this is perfectly fine
and this defines a twisted bimodule structure
$\mathcal{M}=(M^\bullet, \mathbb{M})$, except that $M^\bullet$ is
not finitely generated over $\A(Y)$.

Let us see what the functor $\mathcal{M}_*$ does on
$\Pc_{\A^{0,\bullet}(X\times Y)}$ If $E$ is the smooth sections of
a holomorphic vector bundle, equipped with its
$\mathbb{E}=\dbar$-operator, then $\mathcal{M}_*(E)$ is the
twisted $\A^{0,\bullet}(Y)$-module $E$ (considered only as a
$\A(Y)$-module). The $\mathbb{Z}$-connection
$\mathbb{M}\#\mathbb{E}=(\mathbb{M}\#\mathbb{E})^0 +
(\mathbb{M}\#\mathbb{E})^1$ where $(\mathbb{M}\#\mathbb{E})^0$ is
the  $\dbar$-operator along the fiber and
$(\mathbb{M}\#\mathbb{E})^1$ is the $\dbar$-operator along the
base. If the cohomology of $E$ with respect to
$(\mathbb{M}\#\mathbb{E})^0$ were of locally constant dimension,
then it would define a complex of vector bundles on $Y$ and
$(\mathbb{M}\#\mathbb{E})^1$ would provide a
$\mathbb{Z}$-connection over $Y$ and we would have an object in
$\Pc_{A^{0,\bullet}(Y)}$. In general we would have to resolve
after we take the fiberwise cohomology or perturb before we take
the fiberwise cohomology so that we get a complex of vector
bundles. A pastiori, Grauert's direct image theorem tells us that
we will get something coherent. In this paper, we will merely
leave the answer as the ``quasi"-perfect object $\mathcal{M}_*(E)$
and deal with when the answer lies in the smaller perfect category
in part III.

While the bimodule $M^\bullet$ is not finitely generated over
$\A(Y)$ it is still projective the sense of homological algebra
over topological algebras. We are thus in a relative homological
situation.

We  now set up the general framework. We will work with Fr\'echet
algebras to make things simpler. The correct generalization beyond
Fr\'echet algebras will be bornological algebras as in \cite{Hou}
and has been used recently by \cite{Mey} and \cite{KS}. The
generalization to this context is straightforward. Fix a ground
field $\kk$, either $\mathbb{R}$ or $\mathbb{C}$. A Fr\'echet
space is a locally convex topological vector space over $\kk$
which is defined by a metric invariant under vector space addition
and which is complete. Equivalently, it is a complete locally
convex topological vector space which is defined by a countable
collection of semi-norms $p\in \Lambda$. A Fr\'echet algebra is a
Fr\'echet space $\A$ which is endowed with a $\mathbb{C}$-bilinear
map
\[
\A\times\A\to\A
\]
which is continuous. A Fr\'echet algebra is multiplicatively
convex if for every continuous semi-norm $p\in\Lambda$ we have \[
p(ab)\le p(a)p(b)
\]
This notion becomes important when using perturbation techniques
as in the proof of Grauert's theorem.  A module over a Fr\'echet
algebra $\A$ is a module (for us, usually a right module) $M$
which is a Fr\'echet space such that the module action is
continuous. For two modules $M$ and $N$ over $\A$, a morphism
$T:M\to N$ is a continuous $\A$-linear homomorphism. The space
$\mycal{L}_\A(M,N)$ is a complete locally convex topological
vector space, though no longer Fr\'echet in general. We will
usually consider $\mycal{L}_\A(M,N)$ just as an abstract vector
space. A complex of $\A$-modules is a complex of Fr\'echet spaces
$(M^\bullet, d)$ such that $d$ is continuous and $\A$-linear. If
$(M^\bullet, d)$ and $(N^\bullet,d)$ are complexes of $\A$-modules
then
\[
\mycal{L}_\A^\bullet(M^\bullet,N^\bullet)
\]
is a complex with differential
\[(d\phi)(m)=d(\phi(m))-(-1)^{|\phi|}\phi(d(m)).\]
From now on, when considering Fr\'echet modules we will write
\[\Hom_\A=\mycal{L}_\A  \mbox{  and  } \Hom^\bullet_\A=\mycal{L}^\bullet_\A.\] This is
consistent with our previous notation since for finitely generated
projective modules every $\A$-linear homomorphism is continuous.

We will use $\otimes$ for the completed projective tensor product
of Fr\'echet spaces. (This is usually denoted $\hat{\otimes}$.)
Again, this is consistent with our previous usage since for
finitely generated projective modules, the algebraic tensor
product and the topological tensor product agree.

For a right $\A$-module $M$ and a left $\A$-module $N$, we define
$M\otimes_\A N$ to be the quotient by the closure of the image of
the map
\begin{equation}\label{tensorproduct} M\otimes A\otimes
N\to M\otimes N \,\,\,\,\,\,\,\,\,\,\mbox{ defined by }
\,\,\,\,\,\,\,\, m\otimes a\otimes n\mapsto ma\otimes n-m\otimes
an
\end{equation}
This has the universal property that
$\Hom_\kk(M\otimes_\A N,L)$ is naturally isomorphic to the space
of continuous $\kk$-bilinear maps $T:M\times N\to L$ such that
$T(ma,n)=T(m,an)$. One often defines the tensor product over $\A$
without taking closures. In our case this won't change anything as
will be explained below.

\begin{prop}
Given a Fr\'echet $\A$-module of the form $V\otimes A$, where $V$
is a Fr\'echet space, and any other Fr\'echet module $M$ there is
a canonical isomorphism
\[
\Hom_\A(V\otimes_\kk A, M)\cong \Hom_\kk(V, M)
\]
\end{prop}

An $\A$-module of the form $V\otimes_\kk A$ is called {\em
relatively free}. A module $P$ is {\em relatively projective} if
it is a direct summand  of a relatively free module. That is,
there is a surjection
\[
V\otimes_\kk \A \to P
\]
that admits an $\A$-linear continuous section. From now one, we
use simply free (resp. projective) instead of relatively free
(resp. projective).

The following is standard.
\begin{prop}
Let $P$ be a projective right $\A$-module and $E$ any Fr\'echet
left $\A$-module. Then $P\otimes_\A E$ is equal to what one would
get by taking the quotient in \eqref{tensorproduct} without taking
the closure of the image.
\end{prop}

The category of bounded complexes of projective $\A$-modules and
continuous maps forms a DG-category, $\mathcal{C}_\A$. For
$\A=\kk$ is just the category of bounded complexes of Fr\'echet
spaces with continuous maps. As per usual terminology, a module
over the category $\mathcal{C}_\A$ is a contravariant functor,
that is a functor from the opposite category
$M:\mathcal{C}_\A^\circ\to \mathcal{C}_\kk$. The category of
modules over $\mathcal{C}_\A$ itself forms  a DG-category,
\cite{BK}, $\mathcal{C}_\A^\circ DG-\mbox{Mod}$. The category
$\mathcal{C}_\A$ embeds in $\mathcal{C}_\A^\circ DG-\mbox{Mod}$
via the Yoneda embedding
\[
M^\bullet\mapsto \Hom^\bullet(\cdot, M^\bullet)
\]
Even for bounded complexes of Fr\'echet modules $M^\bullet$ which
are not projective, we can define an object in
$\mathcal{C}_\A^\circ DG-\mbox{Mod}$ using the same formula as the
Yoneda embedding. We note the following obvious proposition that
relates the notion of quasi-isomorphism in our DG-category
$\mathcal{C}_\A^\circ DG-\mbox{Mod}$ to relative homological
algebra.
\begin{prop}
Let
\[
0\to E_1\to E_2\to E_3\to 0
\]
be an exact sequences of not necessarily projective $\A$-modules
and consider it as a complex in $\mathcal{C}_\A^\circ
DG-\mbox{Mod}$. Then it is quasi-isomorphic to $0$ if and only if
it is split as a sequence of topological vector spaces.
\end{prop}

Finally we come to the definition of the quasi-perfect category.
\begin{defn}
Let $(\A^\bullet,d,c)$ (or just $\A^\bullet$ for short) be a
curved DGA, which is Fr\'echet as an algebra and such that $d$ is
continuous in the Fr\'echet topology. We associate to $\A^\bullet$
the DG-category $q\Pc_{\A^\bullet}$  of quasi-perfect twisted
complexes whose objects are $E=(E^\bullet, \mathbb{E})$ where
$E^\bullet$ is a bounded $\mathbb{Z}$-graded right Fr\'echet
$\A$-module which is projective (but not necessarily finitely
generated) and
\[ \mathbb{E}:E^\bullet\to E^\bullet\otimes_\A A^\bullet
\]
is a $\mathbb{Z}$-connection, which is continuous in the
respective Fr\'echet topologies, and it satisfies
\[
\mathbb{E}\circ\mathbb{E}(e)=-e\cdot c
\]

\noindent The morphisms between two such objects
$E_1=(E_1^\bullet, \Eb_1)$ and $E_2=(E_2^\bullet,\Eb_2)$ of degree
$k$ are
\[
\Hom^k_{q\Pc_\A^\bullet}(E_1,E_2)=\{\phi:E_1^\bullet\otimes_\A
\A^\bullet\to E_2^{\bullet}\otimes_\A \A^\bullet\,\,|\,\,
\phi(ea)=(-1)^{k|a|}\phi(e)a\}
\]
where now the morphisms are required to be continuous. The
morphisms are equipped  with a differential defined by
\[
d(\phi)(e)= \Eb_2(\phi(e))-(-1)^{\abs{\phi}}\phi(\Eb_1(e))
\]
When considering cohomology of the $\Hom$-complex, we always
consider the ordinary cohomology of complexes of vector spaces.
That is, we do not quotient out by the closure of the boundaries
and we do not consider any topology on the cohomology spaces. Of
course, the topology of the modules enters when defining what
$\Hom$ is.

 Clearly, $\Pc_{\A^\bullet}$ is a full subcategory of
$q\Pc_{\A^\bullet}$ since a homomorphisms between finitely
generated $\A$-modules is automatically continuous.
\end{defn}
\begin{prop} For a curved DGA $(\A^\bullet, d, c)$ the category $q\Pc_{\A^\bullet}$ is a DG-category.
\end{prop}

We also have the following generalization of a perfect twisted
bimodule. Let $(\A_1^\bullet,d,c_1)$ and $(\A^\bullet_2,d,c_2)$ be
two curved DGAs.
\begin{defn}
A quasi-perfect twisted bimodule is the data
$X=(X^\bullet,\mathbb{X})$ where
\begin{enumerate}
\item $X^\bullet$ is a bounded $\mathbb{Z}$-graded projective
right Fr\'echet $\A_2$-module, (not necessarily finitely
generated) \item $\mathbb{X}:X^\bullet\to
X^\bullet\otimes_{\A_2}\A_2^\bullet$ is a continuous
$\Zb$-connection, \item $\A^\bullet_1$ acts continuously on the
left of $X^\bullet\otimes_{\A_2}\A_2^\bullet$ satisfying

\[
a\cdot (x\cdot b)=(a\cdot x)\cdot b
\]
and
\[
\mathbb{X}(a\cdot(x\otimes b))=da\cdot(x\otimes
b)+a\cdot\mathbb{X}(x\otimes b)
\]
for $a\in \A_1^\bullet$, $x\in X^\bullet$ and $b\in \A_2^\bullet$,
 \item
$\mathbb{X}$ satisfies the following condition
\[
\mathbb{X}\circ\mathbb{X}(x\otimes b)=c_1\cdot (x\otimes
b)-(x\otimes b)\cdot c_2
\]
 on the complex
$X^\bullet\otimes_{\A_2}\A_2^\bullet$. \end{enumerate}\end{defn}

Since tensoring a relatively projective right $A_1$-module by a
$\A_1$-$\A_2$ bimodule which is a relatively project right
$\A_2$-module yields a relatively projective right $\A_2$-module,
we have, as before, that a quasi-perfect twisted bimodule
$\mathcal{X}=(X^\bullet,\mathbb{X})$ defines a DG-functor
\[
\mathcal{X}_*:q\Pc_{\A_1^\bullet}\to q\Pc_{\A_2^\bullet}
\]

\begin{prop}
For two complex manifolds $X$ and $Y$ and $p_2:X\times Y\to Y$ the
projection map, the pair $\mathcal{M}=(M^\bullet, \mathbb{M})$
defined at the beginning of this section defines a
$\A^{0,\bullet}(X\times Y)$-$\A^{0,\bullet}(Y)$ quasi-perfect
twisted bimodule, and thus defines a DG functor
\[
\mathcal{M}_*:q\Pc_{\A^{0,\bullet}(X\times Y) }\to
q\Pc_{\A^{0,\bullet}(Y) }
\]
\end{prop}
\begin{rem}
\begin{enumerate}
\item For an arbitrary holomorphic map $f:X\to Y$ it is possible
to define the pushforward using the techniques at hand. Namely, by
factoring the map as an embedding and a projection, we can resolve
the graph of $f$, and compose with the quasi-perfect twisted
bimodule for the projection. We leave the details as an exercise.
\item The pushforward defined by $\mathcal{M}$ is the derived
pushforward since we always stay in the category $q\Pc$, which
(for the case of the Dolbeault algebra) corresponds to fine
sheaves of projectives.
\end{enumerate}
\end{rem}

By enlarging the category from $\Pc_{\A^\bullet}$ to
$q\Pc_{\A^\bullet}$ we have potentially changed (made more strict)
the notion of quasi-isomorphism between to objects in
$q\Pc_{\A^\bullet}$, even if they are both in $\Pc_{\A^\bullet}$.
However, in the case of a complex manifold, the notion of
quasi-isomorphism between complexes of coherent sheaves is still
preserved in the quasi-perfect category because of the following
result which is an extension of the result for the perfect
category.
\begin{prop}\label{prop:qi_criterion}
Let $E_1$ and $E_2$ be two quasi-perfect twisted complexes. If for
each $p$, $\A^p$ is flat as an $\A$-module, a closed morphism
$\phi\in \Hom^0_{q\Pc_{\A^\bullet}}(E_1,E_2)$ is a
quasi-isomorphism if $\phi^0:(E_1^\bullet,\Eb_1^0)\to
(E_2^\bullet,\Eb_2^0)$ is a quasi-isomorphism of complexes of
$\A$-modules.
\end{prop}

\subsection{DG-equivalences}
We would like to establish a criteria for when two quasi-perfect
twisted bimodules determine a DG-equivalence. Let $(\A^\bullet, d,
c_A)$ and $(\B^\bullet,d, c_B)$ be two curved DGAs. Let us note
that in the curved DGA case, the module $(\A, d)$ is not a perfect
twisted complex over $\A^\bullet$, since it does not satisfy the
correct curvature conditions. It is however, a perfect twisted
bimodule and it is clear that the functors it induces on
$\Pc_{\A^\bullet}$ and $q\Pc_{\A^\bullet}$ are the identity
functors.

Let $(\Pc^\bullet,\mathbb{P})$ be a $\B^\bullet$-$\A^\bullet$
quasi-perfect twisted bimodule and $(\Qc^\bullet,\mathbb{Q})$ an
$\A^\bullet$-$\B^\bullet$ quasi-perfect twisted bimodule. Let us
write $\Xc^\bullet=\Qc^\bullet\otimes_\B\Pc^\bullet$ and
$\Yc^\bullet=\Pc^\bullet\otimes_\A\Qc^\bullet$. Furthermore, let
$\mathbb{X}=\mathbb{Q}\#\mathbb{P}$ and
$\mathbb{Y}=\mathbb{P}\#\mathbb{Q}$. Suppose we have two maps
\[
\alpha: \Xc^\bullet\otimes_\A \A^\bullet\to \A^\bullet
\]
and
\[
\beta: \Yc^\bullet\otimes_\B^\bullet\to \B^\bullet
\]
such that $\alpha$ is a surjective map of $\A^\bullet$-bimodules
and $\beta$ is a surjective map of $\B^\bullet$-bimodules  and
both $\alpha$ and $\beta$ intertwine the $\mathbb{Z}$-connections:
\[
\alpha(\mathbb{X}(x))=d(\alpha(x))
\]
and
\[
\beta(\mathbb{Y}(y))=d(\beta(y))
\]
Under these circumstances, there are natural transformations of
functors
\[
\alpha:\Xc_*\to \mathbb{1}_{q\Pc_{\A^\bullet}} \] and
\[
\beta:\Yc_*\to \mathbb{1}_{q\Pc_{\B^\bullet}}
\]
defined, for example, for $E=(E^\bullet,\mathbb{E})\in
q\Pc_{\A^\bullet}$
\[
\alpha_E:E^\bullet\otimes_\A\Xc^\bullet\otimes_\A\A^\bullet\to
E^\bullet\otimes_\A\A^\bullet
\]
by
\[
\alpha_E(e\otimes x)=e\otimes\alpha(x)
\]
Now according to \eqref{prop:qi_criterion}, to show that
$\Pc^\bullet$ and $\Qc^\bullet$ induce DG-quasi-equivalences, we
need to see that for each $E=(E^\bullet,\mathbb{E})$ that we have
an isomorphism
\[
\alpha_E^0:H^*(E^\bullet\otimes_\A\Xc^\bullet,(\mathbb{E}\#\mathbb{X})^0)\to
H^*(E^\bullet,\mathbb{E}^0)
\]
and similarly for $\Yc_*$. We would like to have a condition that
we can check about $\Xc^\bullet$ and $\Yc_\bullet$ themselves. We
have the following criterion.

Note first, that in the case when the curvature is not zero, that
$\mathbb{X}$ is not a differential. It is not even true that
$(\mathbb{X}^0)^2$ is necessarily zero. One has
$(\mathbb{X}^0)^2=c_A$. Suppose there is an endomorphism
$\Phi:\Xc^\bullet\to \Xc^\bullet$ of degree one, which is
$\A$-linear on the right and the left. Suppose further, that
$\Phi$ satisfies
\begin{equation}\label{eq:gaugechange}
[\mathbb{X}^0,\Phi]+\Phi\circ\Phi=-c_A
\end{equation}
Then of course we can form $\mathbb{X}^0+\Phi$ which now has
square zero and our criterion is
\begin{lem}\label{lem:dgeq}
If $\alpha^0$ induces an isomorphism
\[
H^*(\Xc^\bullet,\mathbb{X}^0+\Phi)\mto{\alpha^0} \A
\]
then the natural transformation $\alpha:\Xc_*\to
\mathbb{1}_{\Pc_{\A^\bullet}}$ is a DG-isomorphism of functors.
\end{lem}
\begin{rem}
As we will see below, one way to interpret the endomorphism $\Phi$
is as a trivialization of a ``gerbe" on the product $(A^\bullet,
d, c_A)\otimes (B^\bullet,d,c_B)$.
\end{rem}

\section{Mukai duality for non-commutative tori}
In this section we state and prove the duality between a
noncommutative complex torus and a gerby complex torus.

\subsection{The Complex noncommutative torus}
We start by describing the  noncommutative tori. We will
describe them in terms of twisted group algebras. Let $V$ be a
real vector space, and $\Lambda\subset V$ a lattice subgroup. The
we can form the group ring $\Sc^*(\Lambda)$, the Schwartz space of
complex valued functions on $\Lambda$ which decrease faster than
any polynomial. Let $B\in \Lambda^2 V^\vee$, and form the
biadditive, antisymmetric group cocycle
$\sigma:\Lambda\times\Lambda\to U(1)$ by \[
\sigma(\lambda_1,\lambda_2)=e^{2\pi i B(\lambda_1,\lambda_2)}
\]
In our computations, we will often implicitly make use of the fact
that $\sigma$ is biadditive and anti-symmetric.  Now we can form
the twisted group algebra $\A(\Lambda;\sigma)$ consisting of the
same space of functions as $\Sc^*(\Lambda)$ but where the
multiplication is defined by
\[
[\lambda_1]\circ[\lambda_2]=\sigma(\lambda_1,\lambda_2)[\lambda_1+\lambda_2]
\]
This is a $*$-algebra where
$f^*(\lambda)=\overline{f(\lambda^{-1})}$. This is one of the
standard ways to describe the (smooth version) of the
noncommutative torus.  Given $\xi\in V^\vee$ it is easy to check
that
\[
\xi(f)(\lambda)=2\pi \imath\langle \xi, \lambda\rangle f(\lambda)
\]
defines a derivation on $\A(\Lambda;\sigma)$. Note that the
derivation $\xi$ is ``real" in the sense that $\xi(f^*)=-\xi(f)$.
Finally define a (de Rham) DGA
\[
\A^\bullet(\Lambda;\sigma)=\A(\Lambda;\sigma)\otimes
\Lambda^\bullet V\] where the differential $d$ is defined on
functions $\phi\in \A(\Lambda;\sigma)$ by
\[
\langle df, \xi\rangle=\xi(f)
\]
for $\xi\in V^\vee$. In other words, for $\lambda\in \lambda$ one
has $d\lambda=2\pi\imath\lambda \otimes D(\lambda)$ where
$D(\lambda)$ denotes $\lambda$ as an element of $\Lambda^1 V$.
Extend $d$ to the rest of $\A^\bullet(\Lambda;\sigma)$ by
Leibnitz. Note that $d^2=0$.

We are most interested in the case where our torus has a complex
structure and in defining the analogue of the Dolbeault DGA for a
noncommutative complex torus.  So now let $V$ will be a vector
space with a complex structure $J:V\to V$, $J^2=-\mathbb{1}$. Let
$g$ be the complex dimension of $V$. Set
$V_\mathbb{C}=V\otimes_\mathbb{R}\mathbb{C}$. Then $J\otimes
1:V_\mathbb{C}\to V_\mathbb{C}$ still squares to $-\mathbb{1}$ and
so $V_\mathbb{C}$ decomposes into $\imath$ and $-\imath$
eigenspaces, $V_{1,0}\oplus V_{0,1}$. The dual $V_\mathbb{C}^\vee$
also decomposes as $V_\mathbb{C}^\vee=V^{1,0}\oplus V^{0,1}$.  Let
$D':V_\mathbb{C}\otimes \mathbb{C}\to V^{1,0}$ and $D'':V\otimes
\mathbb{C}\to V^{0,1}$ denote the corresponding projections.
Explicitly
\[D'=\frac{J\otimes 1 + 1\otimes \imath}{2\imath}\] and
\[D''=\frac{-J\otimes 1 + 1\otimes \imath}{2\imath}\]
and $D=D'+D''$ where $D$ denotes the identity. This also
established a decomposition
\[\Lambda^kV_\mathbb{C}=\otimes_{p+q=\bullet}\Lambda^{p,q}V\]
where $\Lambda^{p,q}V=\Lambda^p V_{1,0}\otimes \Lambda^qV_{0,1}$.
Complex conjugation on $V_\mathbb{C}$ defines an involution and
identifies $V$ with the $v\in V_\mathbb{C}$ such that
$\overline{v}=v$.

Now let $X=V/\Lambda$, a complex torus of dimension $g$,  and
$X^\bullet=\overline{V}^\vee/\Lambda^\vee$ its dual torus. Let
$B\in\Lambda^2 V^\vee$ be a real (constant) two form on $X$. Then
$B$ will decompose in to parts \[ B=B^{2,0}+B^{1,1}+B^{0,2}\]
where $B^{p,q}\in \Lambda^{p,q}V$, $B^{0,2}=\overline{B^{2,0}}$
and $\overline{B^{1,1}}=B^{1,1}$. Now $B^{0,2}\in \Lambda^2
V^{0,1}\cong H^{0,2}(X)$. Then it also represents a class in
$\Pi\in \Lambda^2 V^{0,1}\cong H^0(X^\vee;\Lambda^2 T_{1,0}X)$.
Let $\sigma:\Lambda \wedge \Lambda\to U(1)$ denote the group
$2$-cocycle given by
\[\sigma(\lambda_1,\lambda_2)=e^{2\pi
\imath B(\lambda_1,\lambda_2)}.\] and form as above
$\A(\Lambda;\sigma)$ the twisted group algebra based on rapidly
decreasing functions. Define the Dolbeault DGA
$\A^{0,\bullet}(\Lambda;\sigma)$ to be
\[
\A(\Lambda;\sigma)\otimes \Lambda^\bullet V_{1,0}
\]
where for $\lambda\in \A(\Lambda;\sigma)$ we define
\[
\dbar\lambda=2\pi\imath \lambda\otimes D'(\lambda)\in
\A(\Lambda;\sigma)\otimes V_{1,0}
\]
We can then extend $\dbar$ to the rest of
$\A^{0,\bullet}(\Lambda;\sigma) $ by the Leibnitz rule.

\begin{rem} Let us make a comment about this definition. Even though we
are defining the $\dbar$ operator, we are using the $(1,0)$
component of $V_\mathbb{C}$. This is because of the duality. In
the case of the trivial cocycle $\sigma$, this definition is meant
to reconstruct the Dolbeault algebra on $X^\vee$. In this case,
$\A^{0,\bullet}(X^\vee)\cong \A(X^\vee)\otimes \Lambda^\bullet
T^{0,1}_0X^\vee$. But
\[T^{0,1}_0X^\vee = \overline{\overline{V}^\vee}^\vee\cong
V_{1,0}.\] \end{rem}
To check the reasonableness of this definition
we have
\begin{prop} If $\sigma=1$ is the trivial
cocycle, then the DGA $(\A^{0,\bullet}(\Lambda;\sigma), \dbar)$ is
isomorphic to the Dolbeault DGA $(\A^{0,\bullet}(X^\vee),\dbar)$.
\end{prop}

\subsection{The holomorphic gerby torus}
The other side of the duality involves a holomorphic
$\mathcal{O}^\times$ gerbe over the dual of $X^\vee$, that is,
over $X$. A gerbe is usually described in terms of a $2$-cocycle
with values in $\mathcal{O}^\times$. Our gerbes are topologically
trivial, though not holomorphically so. Using this topological
trivializtion, we can express it in terms of a very simple curved
DGA. See \cite{BD} for the precise relationship between a
holomorphic gerbe defined in terms of a cocycle and the curved DGA
we now describe.

Consider $B\in \Lambda^2 V^\vee$ decomposed as
$B=B^{0,2}+B^{1,1}+B^{2,0}$ as above. Let $\A^{0,\bullet}(X;B)$
denote the curved DGA $(\A^{0,\bullet}(X),\dbar,2\pi\imath
B^{0,2})$. (We write $A^{0,\bullet}(X;B)$ even though we have only
used the $(0,2)$ component. If we had constructed the bigraded
Dolbeault algebra, $\A^{\bullet,\bullet}(X;B)$, we would have used
all the components of $B$.) While the underlying DGA (non-curved)
of $\A^{0,\bullet}(X,B)$ is the same as the Dolbeault algebra of
$X$, the existence of the curvature changes the notion of a module
and thus the category $\Pc$ and $q\Pc$.

We will show there is a DG-quasi-equivalence of categories between
$q\Pc_{\A^{0,\bullet}(\Lambda;\sigma)}$ and
$q\Pc_{\A^{0,\bullet}(X;B)}$. This will be implemented by a pair
of quasi-perfect twisted  bimodules which are deformed versions of
the Poincare sheaves.

\subsection{The deformed Poincare line bundles}
A reader familiar with $C^*$-algebra $K$-theoretic techniques will
notice the similarity of the constructions in this section with
the Kasparov bimodules that implement the Baum-Connes assembly
map. On the other hand, our constructions are also deformations of
the Poincare sheaves on $X\times X^\vee$ realized as a line bundle
with a $\dbar$-connection.

We define an $\A(X)$-$\A(\Lambda;\sigma)$ bimodule $\Pc$ by
setting $\Pc$ to be the vector space of Schwartz functions
$\Sc(V)$ with the left action of $\A(X)$ to be just pulling back a
function from $X$ to $V$ and multiplying. The right action of
$\A(\Lambda;\sigma)$ is defined for $\lambda\in
\Lambda\subset\A(\Lambda;\sigma)$ and for $p\in\Pc$ by
\[
p\cdot\lambda(v)= \sigma(\lambda,v)p(\lambda+v)
\]
Here $\sigma$ has been extended to map from $\Lambda\times V\to
U(1)$ using the same formula: $\sigma(\lambda,v)=e^{2\pi\imath
B(\lambda,v)}.$ This extension still satisfies the obvious
``cocycle" relation
\[
\delta\sigma(\lambda_1,\lambda_2,v)=1
\]
where
\[
\delta\sigma(\lambda_1,\lambda_2,v)=\sigma(\lambda_2,v)\sigma(\lambda_1+\lambda_2,v)^{-1}
\sigma(\lambda_1,\lambda_2+v)\sigma(\lambda_1,\lambda_2)^{-1}
\] One checks easily that this makes $\Pc$ into an
$\A(X)$-$\A(\Lambda;\sigma)$ bimodule.

\begin{lem}\label{proj-lemma}
The module $\Pc$ is projective as a right $\A(\Lambda;\sigma)$
module. It is also projective as a left $\A(X)$-module. \end{lem}
\proof The proof of this is rather standard and uses a ``partition
of unity" $h:V\to\mathbb{R}$ a compactly supported nonnegative
$C^\infty$ function such that
\[
 \sum_{\gamma\in \Lambda}h(\gamma + v)=1
\]
for all $v\in V$. Then we use it to split the map $\Pc\otimes
\A(\Lambda;\sigma)\to \Pc$ given by the action. Define the
$\A(\Lambda;\sigma)$-module map $\iota:\Pc\to \Pc\otimes
\A(\Lambda;\sigma)$ (a free $\A(\Lambda;\sigma)$ module)
\begin{equation}
\iota p=\sum_{g\in\Lambda} (p\cdot (-g))(v)h(v)\otimes g
\end{equation}
It is easy to check this is a splitting. We check,  for $p\in\Pc$
and $\mu\in \Lambda$ that
\begin{equation}\label{proj-1}\begin{split}
\iota(p\cdot\mu)(v) & =\sum_{g\in\Lambda} ((p\cdot\mu)\cdot
(-g))(v)h(v)\otimes g \\
 & = \sum_{g\in\Lambda} \sigma(\mu,-g)(p\cdot(\mu -g)
)(v)h(v)\otimes g \end{split}\end{equation} using the fact that
$(p \cdot\mu)\cdot g)=\sigma(\mu,g)(p\cdot(\mu + g))$. Letting
$-\tau=\mu -g$ so $g=\tau+\mu$ \eqref{proj-1} becomes
\begin{equation}\label{proj-2}
=\sum_{\tau\in\Lambda}  (p\cdot
(-\tau))(v)h(v)\otimes \sigma(\mu+\tau,\mu)(\tau+\mu) \\
= (\iota(p)\cdot\mu)(v)
\end{equation}
Thus $\iota$ is a module homomorphism.

The projectivity as a left $\A(X)$-module is even easier. $\Pc$ is
the global $C^\infty$ sections of an infinite dimensional
Fr\'echet space bundle over $X$, \[ \Pc\cong \Gamma(X;
\mathcal{V})
\] where the total space of this bundle is
\[
\mathcal{V}=V\times_\Lambda \mycal{S}^*\Lambda
\]
Then the standard fact that sections of a bundle (albeit infinite
dimensional) is projective still holds.
 \qed

 To define our quasi-perfect twisted  bimodule we set
\[
\Pc^\bullet=\Pc\otimes \Lambda^\bullet V^{0,1}
\]
The actions extend in an obvious way. The projectivity follows
from Lemma \ref{proj-lemma}.

$B^{0,2}$ defines a $\dbar$-closed $(0,2)$ form on $X$. Its
pullback to $V$ is $\dbar$ exact. In fact, $B^{0,2}=\dbar \omega$
where $\omega$ is a $(0,1)$ form on $V$ which can be described as
follows. In coordinates, using the reality of $B$ we have (always
use the summation convention)
\[ B=
b_{ij}dz_idz_j+\overline{b_{ij}}d\overline{z}_id\overline{z}_j+c_{ij}dz_id\overline{z}_j\]
where $c_{ij}$ is a skew Hermitian matrix, and we may (and do)
assume that $b_{ij}$ is skew symmetric. Set
\begin{equation}\label{omega}
\omega=\overline{b_{ij}z_i}d\overline{z}_j+c_{ij}z_id\overline{z}_j
\end{equation}
Then $\dbar\omega=B^{0,2}$. We now observe the following
relationships between $\omega$ and $\sigma$. Let $\lambda\in
\Lambda$. Then

\begin{equation}\begin{split}
\dbar(\sigma(\lambda,\cdot))& =\dbar (e^{2\pi\imath
(b_{ij}\lambda_iz_j+\overline{b_{ij}
\lambda_iz_j}+c_{ij}(\lambda_i\overline{z}_j-
z_i\overline{\lambda}_j)})\\
      &=2\pi\imath(\overline{b}_{ij}\overline{\lambda}_id\overline{z}_j
      +c_{ij}\lambda_id\overline{z}_j)\sigma(\lambda,z)\end{split}\end{equation}

\noindent  On the other hand, one easily calculates
\begin{equation}\label{omega1}
\begin{split}
2\pi\imath(\omega-r_\lambda^*\omega)\sigma(\lambda,\cdot) & =
2\pi\imath(\overline{b}_{ij}\overline{z}_id\overline{z}_j+c_{ij}z_id\overline{z}_j
-\overline{b}_{ij}\overline{(z_i+\lambda_i)}d\overline{z}_j-c_{ij}(z_i+\lambda_i)d\overline{z}_j)\sigma(\lambda,\cdot)\\
    & =
-\dbar\sigma(\lambda,\cdot)\end{split}\end{equation}
 where
$r_\lambda(v)=v+\lambda$. Similarly, letting
$l_\lambda(v)=v-\lambda$ we have
\begin{equation}\label{omega2}
2\pi\imath(\omega-l_\lambda^*\omega)\sigma(\lambda,\cdot)
=\dbar\sigma(\lambda,\cdot).\end{equation}

Define a $\Zb$-connection $\mathbb{P}$ on $\Pc^\bullet$,
\[\mathbb{P}:\Pc^\bullet\to \Pc^\bullet\otimes_{\A(\Lambda;\sigma)}\A^{0,\bullet}(\Lambda;\sigma)\cong
\Pc^\bullet\otimes_\mathbb{C}\Lambda^\bullet V_{1,0}\] as
$\mathbb{P}=\mathbb{P}^0+\mathbb{P}^1$ where
\[
\mathbb{P}^0(p)(v)=\dbar_V(p)(v)+2\pi\imath \omega(v)\wedge p(v)
\]
and
\[
\mathbb{P}^1p(v)=-2\pi\imath p(v)D'(v)
\]
Let us check that this is indeed a $\Zb$-connection. There are two
things to check.

First we check that $\mathbb{P}^0$ is $\A(\Lambda;\sigma)$ linear:
For $\lambda\in\Lambda\subset\A(\Lambda;\sigma)$ and $p\in
\Pc^\bullet$ we have
\begin{equation}
\begin{split}
\mathbb{P}^0(p\cdot\lambda)=&\dbar_V(p\cdot\lambda)+2\pi\imath
\omega\wedge (p\cdot\lambda) \\
 = & \dbar_V(r_\lambda^*p\cdot
 \sigma(\lambda,\cdot))+2\pi\imath\omega\wedge r_\lambda^*p\cdot
 \sigma(\lambda,\cdot)\\
=  &
\dbar_V(r_\lambda^*p)\sigma(\lambda,\cdot)+r_\lambda^*p\dbar_V(\sigma(\lambda,\cdot))+2\pi\imath
\omega\wedge r_\lambda^*p\cdot \sigma(\lambda,\cdot) \\
=  &
r_\lambda^*\dbar_Vp\cdot\sigma(\lambda,\cdot)+r_\lambda^*p\wedge
2\pi\imath(r_\lambda^*\omega-\omega)\sigma(\lambda,\cdot)+2\pi\imath\omega\wedge
r_\lambda^*p \cdot\sigma(\lambda,\cdot)\\
= & r_\lambda^*(\dbar_V(p)+2\pi\imath\omega\wedge
p)\sigma(\lambda,\cdot) \\
 = & \mathbb{P}^0(p)\cdot \lambda
\end{split}
\end{equation}

Second we need that   $\mathbb{P}^1$ satisfies Leibnitz with
respect to $\A(\Lambda;\sigma)$: Again, for
$\lambda\in\Lambda\subset\A(\Lambda;\sigma)$ and $p\in
\Pc^\bullet$ we have
\begin{equation}
\begin{split}
\mathbb{P}^1(p\cdot\lambda)(v) =& -2\pi\imath
(p\cdot\lambda)D'(v)\\
 = & -2\pi\imath
r_\lambda^*p(v)D'(v)\sigma(\lambda,v)
\end{split}
\end{equation}
while
\begin{equation}\begin{split}
\mathbb{P}^1(p)\cdot\lambda(v)+p\dbar(\lambda)(v) \\
= & (-2\pi\imath pv) D'(v))\cdot\lambda +2\pi\imath p\cdot \lambda
D'(\lambda) \\
= & -2\pi\imath
r_\lambda^*p(v)D'(v+\lambda)\sigma(\lambda,v)+2\pi\imath
r_\lambda^*p(v)D'(\lambda)\sigma(\lambda,v)\\
 = & -2\pi\imath r_\lambda^*p(v) D'(v)\sigma(\lambda,v)
\end{split}\end{equation}
Now, $\A^{0,\bullet}(X;B)$ acts on the left of
$\Pc^\bullet\otimes\A^{0,\bullet}(\Lambda;\sigma)$ through its
action on $\Pc^\bullet$. The fact that for $\eta\in
\A^{0,\bullet}(X;B)$ and $p\in \Pc^\bullet$ satisfies
\[
\mathbb{P}(\eta\cdot p)=\dbar_X\eta\cdot p+\eta\cdot\mathbb{P}(p)
\]
is easy to verify.

 Finally, let us note that
\[\begin{split}
\mathbb{P}(\mathbb{P}(p))= &
\mathbb{P}(\dbar_Vp+2\pi\imath\omega\wedge p-2\pi\imath p\otimes
D'\\
 = & \dbar_V(\dbar_Vp+2\pi\imath\omega\wedge p-2\pi\imath p\otimes
D') \\
   + &  2\pi\imath\omega\wedge(\dbar_Vp+2\pi\imath\omega\wedge p-2\pi\imath p\otimes
D')\\
 - & 2\pi\imath (\dbar_Vp+2\pi\imath\omega\wedge p-2\pi\imath p\otimes
D')\wedge D' \\
 = & \dbar_V ^2 p +2\pi\imath \dbar_V\omega\wedge p -2\pi\imath
 \omega\wedge \dbar_V p-2\pi\imath \dbar_V p\otimes D'-2\pi\imath
 p\otimes \dbar_V D' \\
 + &  2\pi\imath\omega\wedge(\dbar_Vp+2\pi\imath\omega\wedge p-2\pi\imath p\otimes
D')\\
-  & 2\pi\imath (\dbar_Vp+2\pi\imath\omega\wedge p-2\pi\imath
p\otimes D')\wedge D'
 \end{split}
 \]
 After the obvious cancellations we are left with
 \[
2\pi\imath \,\,\dbar_V\omega\wedge p-2\pi\imath p\otimes \dbar_V
D'
 \]
Now $\dbar_V\omega=B^{0,2}$ and as a map from $V\to \Lambda^1
V_{1,0}$, $D'$ is holomorphic, so $\dbar_V D'=0$. Thus,
\begin{equation}\label{P-curv}
\mathbb{P}(\mathbb{P}p)=2\pi\imath  B^{0,2}\wedge p
\end{equation}
We thus arrive at the
\begin{prop}
$(\Pc^\bullet,\mathbb{P})$ forms a
$\A^{0,\bullet}(X;B)-\A^{0,\bullet}(\Lambda;\sigma)$quasi-perfect
twisted bimodule.
\end{prop}

We define a bi-twisted  complex that implements a DG-functor in
the opposite direction. Set $\Qc^\bullet=\Sc(V;\Lambda^\bullet
V_{1,0})$. In this case we define a left
$\A(\Lambda;\sigma)$-action by defining for
$\lambda\in\Lambda\subset \A(\Lambda;\sigma)$ and $q\in
\Qc^\bullet$ the action
\[
\lambda\cdot q(v)=\sigma(\lambda,-\lambda+v)q(-\lambda+v)
\]
Note that in our case
$\sigma(\lambda,-\lambda+v)=\sigma(\lambda,v)$ but we have written
it as above because it is the correct formula, and it works in
more general situations. It is straightforward to check that this
is indeed a left action. The action extends to the rest of the DGA
$\A^{0,\bullet}(\Lambda;\sigma)$ in the obvious way. $\A(X)$ acts
again by pull pulling a function on $X$ up to $V$ and multiplying.
Exactly as for the case of $\Pc^\bullet$, $\Qc^\bullet$ is also
projective on both sides.

Define a $\mathbb{Z}$-connection
$\mathbb{Q}=\mathbb{Q}^0+\mathbb{Q}^1$ where
\[
(\mathbb{Q}^0q)(v)=2\pi\imath q(v)D'(v)
\]
and
\[
(\mathbb{Q}^1q)(v)=\dbar_Vq- q\wedge 2\pi\imath\omega
\]
Calculations similar to the ones for $\mathbb{P}$ show that
$\mathbb{Q}$ is a $\Zb$-connection and that
\[
\mathbb{Q}(\mathbb{Q}q)=q\wedge(-2\pi\imath B^{0,2})
\]
Hence
\begin{prop}
The pair $(\Qc^\bullet,\mathbb{Q})$ forms a
$\A^{0,\bullet}(\Lambda;\sigma)-\A^{0,\bullet}(X;B)$ quasi-perfect
twisted bimodule.
\end{prop}

Define a homomorphism of $\A(\Lambda;\sigma)-\A(\Lambda;\sigma)$
bimodules
\begin{equation}\label{alpha}
\alpha:\Qc\otimes_{\A(X;B)}\Pc\to \A(\Lambda;\sigma)
\end{equation}
by
\[
\alpha(q\otimes p)= \sum_\lambda \left[\int_V
q(v+\lambda)p(v)\sigma^{-1}(\lambda,v)dv \right][\lambda]
\]
Clearly, $\alpha (q \phi\otimes p)=\alpha (q\otimes  \phi p)$ for
$\phi\in \A(X;B)$.  We check that this is a map of
$\A(\Lambda;\sigma)$ bimodules. We check the compatibility with
the right  action. The left action is slightly simpler.
\begin{equation}\begin{split}
\alpha( q\otimes p\cdot\mu) =& \sum_\lambda \left[\int_V
q(v+\lambda)(p\cdot \mu)(v)\sigma^{-1}(\lambda,v)dv
\right][\lambda] \\
 =& \sum_\lambda \left[\int_V
q(v+\lambda)p(v+\mu)\sigma(\mu,v)\sigma^{-1}(\lambda,v)dv
\right][\lambda]\\
\mbox{Setting } w=v+\mu \mbox{ we get } &  \sum_\lambda
\left[\int_V q(w+\lambda-\mu
)p(w)\sigma(\mu,w-\mu)\sigma^{-1}(\lambda,w-\mu)dw
\right][\lambda]\\
\mbox{Substituting } \tau=\lambda-\mu \mbox{ we get } &\sum_\tau
\left[\int_V
q(w+\tau)p(w)\sigma(\mu-\lambda,w)\sigma^{-1}(\lambda,-\mu)dw
\right][\tau+\mu ]\\
 =& \sum_\tau \left[\int_V
 q(w+\tau)p(w)\sigma^{-1}(\tau,w)dw\right][\tau+\mu]\sigma(\tau,\mu)\\
  =& \alpha( q\otimes p)\cdot \mu \end{split}
  \end{equation}

\subsection{The duality}
Our main theorem is
\begin{thm}\label{theorem-main}
The quasi-perfect twisted bimodules $\mathcal{P}$ and
$\mathcal{Q}$ define DG functors
\[
\mathcal{P}_*:q\Pc_{\A^{0,\bullet}(X;B)}\leftrightarrows
q\Pc_{\A^{0,\bullet}(\Lambda,\sigma)}:\mathcal{Q}_*
\]
which implement DG-quasi-equivalences of DG-categories.
\end{thm}

The rest of this section will be devoted to proving this theorem.
We prove it by calculating the composition of the twisted
bimodules and showing that in each direction they induce functors
equivalent to the identity functor (shifted by $g$).

Write $\mathbb{X}$ for the composed $\Zb$-connection
$\mathbb{Q}\#\mathbb{P}$ on the quasi-perfect twisted bimodule
$\Xc^\bullet=\Qc^\bullet\otimes_{\A(X;B)}\Pc^\bullet$. We
calculate the zero component
$\mathbb{X}^0=(\mathbb{Q}\#\mathbb{P})^0=\mathbb{Q}\otimes
\mathbb{1}+\mathbb{1}\otimes \mathbb{P}^0$.
\newline We have
$(\mathbb{Q}\otimes \mathbb{1}+\mathbb{1}\otimes
\mathbb{P}^0)(q\otimes p)=$
\begin{equation}\label{a1}
 \cf qD'\otimes p +\dbar_V(q)\otimes p -\cf q\wedge \omega
\otimes p +q\otimes \dbar_V(p)+\cf q\otimes \omega\wedge p
\end{equation}
Now $\Xc^\bullet \cong \Sc(V\times \Lambda;\Lambda^\bullet
V_{1,0}\otimes \Lambda^\bullet V^{0,1})$ via the isomorphism from
$\Sc(V)\otimes_{\A(X;B)}\Sc(V)\cong \Sc(V\times_X V)\cong
\Sc(V\times \Lambda)$ sending $(v_1,v_2)\mapsto (v_1,v_2-v_1)$.
Under this isomorphism, the left and right $\A(\Lambda;\sigma)$-
actions can be written as \begin{equation}\begin{split} (\mu\cdot
\phi)(z,\lambda) & =\sigma(\mu,z+\lambda)\phi(z-\mu,\lambda+\mu)
\\
(\phi\cdot \mu)(z,\lambda)&=\sigma(\mu,z)\phi(z,\lambda+\mu)
\end{split}\end{equation}
Furthermore, using this isomorphism and writing \eqref{a1} in
coordinates we can rewrite for $\phi\in \Sc(V\times
\Lambda;\Lambda^\bullet V_{1,0}\otimes \Lambda^\bullet V^{0,1})$

\begin{equation}\label{a2}
\mathbb{X}^0(\phi)(z,\lambda)=\sum_{j=1}^g\left[\cf z_j
d\overline{\zeta_j}\wedge\phi
 +d\overline{z_j}\wedge\frac{\partial}{\partial \overline{z_j}}\phi
 +\cf\sum_i( \overline{\lambda_ib_{ij}}+c_{ij}\lambda_i)d\overline{z_j}\wedge\phi
\right]
\end{equation}
Here we are writing $D':V\to\Lambda V_{1,0}$ in coordinates as
\[
D'(z)=\sum_j z_jd\overline{\zeta_j}
\]
where the notation $d\overline{\zeta_j}$ is a basis for $V_{1,0}$.
We write it this way since they are really anti-holomorphic basis
of $\overline{\overline{V}^\vee}^\vee$. If we write
\begin{equation}\label{def:B_j}
B_j(\lambda)=\sum_i \lambda_ib_{ij}+\overline{\lambda_ic_{ij}}
\end{equation}
then \eqref{a1} can be written
\begin{equation}\label{a2}
\mathbb{X}^0(\phi)(z,\lambda)=\sum_j\cf
z_jd\overline{\zeta_j}\wedge\phi
 +d\overline{z_j}\wedge\frac{\partial}{\partial \overline{z_j}}\phi
 +\cf \overline{B_j(\lambda)}d\overline{z_j}\wedge\phi
\end{equation}

Now we equip $\Xc^\bullet$ with an inner product. First, define on
$\Lambda^\bullet V_{1,0}\otimes \Lambda^\bullet V^{0,1}$ a
Hermitian product by declaring
$d\overline{z}_j,d\overline{\zeta}_k$ to be orthonormal.  Then for
$\phi_1, \phi_2\in \Xc^\bullet$ set
\begin{equation}\label{innerproduct-defn}
\langle \phi_1,\phi_2\rangle=\int_{V\times \Lambda} \langle
\phi_1, \phi_2\rangle dv\end{equation} where $dv$ is Lebesgue
measure on $V$ times counting measure on $\Lambda$. We calculate
the corresponding Laplacian with respect to this inner product,
$(\mathbb{X}^0)^*\mathbb{X}^0+\mathbb{X}^0(\mathbb{X}^0)^*$. The
adjoint \[ (\mathbb{X}^0)^*(\phi)(z,\lambda)=\sum_j -\cf
\overline{z_j}\iota_{\frac{\partial}{\partial
\overline{\zeta_j}}}\phi  -\frac{\partial}{\partial
z_j}\iota_{\frac{\partial}{\partial \overline{z_j}}}\phi
 -\cf
B_j(\lambda)\iota_{\frac{\partial}{\partial \overline{z_j}}}\phi
\]
($\iota_\xi$ denotes contraction with respect to the vector field
$\xi$.)  We first calculate $(\mathbb{X}^0)^*\mathbb{X}^0$ on
functions. It follows from \eqref{a2} that for $\phi$ a function
that
\begin{equation}\label{a3}
(\mathbb{X}^0)^*\mathbb{X}^0\phi (z,\lambda)=
\sum_j\left(-\frac{\partial^2}{\partial z_j\partial
\overline{z_j}}-\cf
\overline{B_j(\lambda)}\frac{\partial}{\partial z_j}-\cf
B_j(\lambda)\frac{\partial}{\partial
\overline{z_j}}+4\pi^2(|z_j|^2+|B_j(\lambda)|^2)\right)\phi
\end{equation}
Write $Y_j(\lambda)$ for the first order differential operator
\[
Y_j(\lambda)(\phi)=\overline{B_j(\lambda)}\frac{\partial}{\partial
z_j}\phi+ B_j(\lambda)\frac{\partial}{\partial \overline{z_j}}\phi
\]
We can then rewrite \eqref{a3} as
\begin{equation}\label{a3.5}
(\mathbb{X}^0)^*\mathbb{X}^0\phi (z,\lambda)=
\sum_j\left(-\frac{\partial^2}{\partial z_j\partial
\overline{z_j}}- \cf
Y_j(\lambda)+4\pi^2(|z_j|^2+|B_j(\lambda)|^2)\right)\phi
\end{equation}
 Note that
$Y_j(\lambda)(z)=\overline{B_j(\lambda)}$ and
$Y_j(\lambda)(\overline{z})=B_j(\lambda).$ Define the deformed
Gaussian, for $\mu \in \Lambda$
\[
b_\mu(z,\lambda)=\left\{\begin{array}{cc}
\exp\left(-2\pi(|z|^2+\sum_j\imath B_j(\lambda)z_j+\sum_j\imath\,
\overline{B_j(\lambda)z_j})\right)& \mbox{ if }\lambda= \mu \\
0 & \mbox{ if } \lambda\ne \mu\end{array}\right. \] We calculate
\begin{equation}\label{a4}
\begin{split}
-\frac{\partial^2}{\partial z_j\partial
\overline{z_j}}b_\mu(z,\lambda) &= -\frac{\partial}{\partial
z_j}(-2\pi b_\mu(z,\mu) (z_j+\imath\overline{B_j(\mu)}))
\\
 & =-4\pi^2b_\mu(z,\mu)(\overline{z_j}+\imath
B_j(\mu))(z_j+\imath\,\overline{B_j(\mu)})+2\pi b_\mu \\
 & = -4\pi^2 b_\mu(z,\mu) (|z_j|^2+\imath B_j(\mu)z_j +\imath\,
 \overline{B_j(\mu)z_j}-|B_j(\mu)|^2)+2\pi b_\mu
 \end{split}\end{equation}
 if $\lambda=\mu$ and is $0$ for $\lambda\ne\mu$. And we see that
 \begin{equation}\label{a5}
 \begin{split}
  -\cf Y_j(\lambda)b_\mu(z,\lambda) & = -\cf (-2\pi b_\mu)Y_j(\lambda)\left(|z_j|^2+\imath
  B_j(\mu)z_j+\imath\,
 \overline{B_j(\mu)z_j} \right)\\
   & = 4\pi^2\imath
   b_\mu(z,\mu)\left(\overline{B_j(\mu)z_j}+B_j(\mu)z_j+\imath
  \overline{B_j(\mu)}B_j(\mu)+\imath
   B_j(\mu)\overline{B_j(\mu)}\right)\\
    & = 4\pi^2b_\mu(z,\mu)\left(\imath \,\overline{B_j(\mu)z_j}+\imath
    B_j(\mu)z_j-2|B_j(\mu)|^2\right)
\end{split}
\end{equation}
if $\lambda=\mu$ and is $0$ for $\lambda\ne\mu$. Adding \eqref{a4}
and \eqref{a5} together we get
\[
-4\pi^2b_\mu(z,\mu)\left (|z_j|^2+|B_j(\mu)|^2\right) +2\pi
b_\mu(z,\mu)
\]
if $\lambda=\mu$ and is $0$ for $\lambda\ne\mu$. And finally we
get that
\begin{equation}\label{a6}
\begin{split}
(\mathbb{X}^0)^*\mathbb{X}^0b_\mu(z,\lambda)& =
\sum_{j=1}^g\left(-\frac{\partial^2}{\partial z_j\partial
\overline{z_j}}- \cf
Y_j(\lambda)+4\pi^2(|z_j|^2+|B_j(\lambda)|^2)\right)b_\mu \\
 & = \sum_j -4\pi^2\left(|z_j|^2+|B_j(\lambda)|^2\right)b_\mu +2\pi
 b_\mu
 +4\pi^2\left(|z_j|^2+|B_j(\lambda)|^2)\right)b_\mu \\
  & = 2\pi g b_\mu
\end{split}\end{equation}

The full Laplacian
$\mathbb{\Delta}^0=(\mathbb{X}^0)^*\mathbb{X}^0+\mathbb{X}^0(\mathbb{X}^0)^*$
can be calculated in a straightforward manner  as

\begin{equation}\label{b1}
\mathbb{\Delta}^0\phi (z,\lambda)
 = \sum_{j=1}^g\left(-\frac{\partial^2}{\partial
z_j\partial \overline{z_j}}- \cf
Y_j(\lambda)+4\pi^2(|z_j|^2+|B_j(\lambda)|^2)-\cf
(d\overline{z}_j\circ
\iota_{\frac{\partial}{\partial\overline{\zeta}_j}}+\iota_{\frac{\partial}
{\partial\overline{z}_j}}\circ d\overline{\zeta}_j)\right) \phi.
\end{equation}
Let us call $L_j=d\overline{z}_j\circ
\iota_{\frac{\partial}{\partial\overline{\zeta}_j}}+\iota_{\frac{\partial}
{\partial\overline{z}_j}}\circ d\overline{\zeta}_j$. Then we can
find a basis of eigenvectors  for $L=\sum_jL_j$ acting on
$\Lambda^\bullet V_{1,0}\otimes \Lambda^\bullet V^{0,1}$. Set
\[
e_j^\pm=d\overline{z}_j\pm \imath d\overline{\zeta}_j
\]
Then
\[
L(e_j^\pm)=\pm\imath e_j^\pm
\]
And more generally, for  $I=(i_1< i_2< \cdots < i_k)$ and $J=(j_1
< j_2< \cdots < j_l)$ we have
\[
L(e^+_I\wedge e^-_J)=(k-l)\imath(e^+_I\wedge e^-_J)
\]
Hence we have the eigenvector decomposition
\[
\Lambda^\bullet V_{1,0}\otimes \Lambda^\bullet
V^{0,1}=\oplus_{I,J}\mbox{ span }e^+_I\wedge e^-_J
\]
and on the $e^+_I\wedge e^-_J$ component, we have
$\mathbb{\Delta}^0\phi (z,\lambda)$
\begin{equation}\label{b1}
= \left\{\sum_j\left(-\frac{\partial^2}{\partial z_j\partial
\overline{z_j}}- \cf
Y_j(\lambda)+4\pi^2(|z_j|^2+|B_j(\lambda)|^2)\right)+2\pi(k-l)\right\}
\phi.
\end{equation}

We now solve this deformed Harmonic oscillator. Indeed the
solution is a deformation of the classical solution of the
Harmonic oscillator by Hermite functions, constructed using
creation and annihilation operators. We will carry out the details
in the case of one complex dimension, the higher dimensional case
following easily because the variables in our equation are all
separable.

Thus, consider the operator
\begin{equation}\label{hamiltonian}
\mathbb{H}\phi=\left(-\frac{\partial^2}{\partial z \partial
\overline{z}}-\cf (B(\lambda)\frac{\partial}{\partial
\overline{z}}+\overline{B(\lambda)}\frac{\partial}{\partial
z})+4\pi^2(|z|^2+|B(\lambda)|^2)\right)\phi
\end{equation}
defined on the Schwarz space $\Sc(\mathbb{C}\times\Lambda)$. We
calculate the eigenvalues and eigenvectors of $\mathbb{H}$. Define
operators
\begin{equation}\label{anncreops}
\begin{split} \mathbb{A}\phi(z,\lambda) & =\left(\frac{\partial}{\partial
\overline{z}}+2\pi(z+\imath \,\,\overline{B(\lambda)})\right)\phi\\
 \mathbb{A}^*\phi(z,\lambda) & =\left(-\frac{\partial}{\partial
z}+2\pi(\overline{z}-\imath B(\lambda))\right)\phi\\
 \mathbb{B}\phi (z,\lambda)& =\left(\frac{\partial}{\partial
z}+2\pi(\overline{z}+\imath B(\lambda))\right)\phi \\
 \mathbb{B}^*\phi(z,\lambda) & =\left(-\frac{\partial}{\partial
\overline{z}}+2\pi(z-\imath \,\,\overline{B(\lambda)})\right)\phi
\end{split}\end{equation}
These operators satisfy the following commutation relations.
\begin{equation}\label{commeq}
\begin{split}
[\mathbb{H},\mathbb{A}]=-2\pi\mathbb{A}, & \,\,\,\,
[\mathbb{H},\mathbb{A}^*]=2\pi\mathbb{A}^*\\
[\mathbb{H},\mathbb{B}]=-2\pi\mathbb{B}, &\,\,\,\,
[\mathbb{H},\mathbb{B}^*]=2\pi\mathbb{B}^*\\
[\mathbb{A},\mathbb{A}^*]=4\pi, &\,\,\,\,
[\mathbb{B},\mathbb{B}^*]=4\pi \\
[\mathbb{A},\mathbb{B}]=0, & \,\,\,\,[\mathbb{A},\mathbb{B}^*]=0\\
[\mathbb{A}^*,\mathbb{B}]=0,  &\,\,\,\,
[\mathbb{A}^*,\mathbb{B}^*]=0
\end{split}\end{equation}
Set $b_\mu^{0,0}(z,\lambda)=b_\mu(z,\lambda)$ from above. Define
recursively,
\begin{equation}\label{Hermite}
\begin{split}
b_\mu^{i+1,j}(z,\lambda)=\mathbb{A}^*b_\mu^{i,j}(z,\lambda), &
\,\,\,\,
b_\mu^{i,j+1}(z,\lambda)=\mathbb{B}^*b_\mu^{i,j}(z,\lambda)
\end{split}
\end{equation}
This is well defined since $\mathbb{A}^*$ and $\mathbb{B}^*$
commutate. Moreover,
$\mathbb{A}b_\mu^{0,0}=\mathbb{B}b_\mu^{0,0}=0$. The following
follows from well-known techniques as in \cite{RS}.
\begin{thm}
The functions $b_\mu^{i,j}\in \Sc(\mathbb{C}\times\Lambda)$ form
an orthogonal complete basis of the closure
$L^2(\mathbb{C}\times\Lambda)$ of $\Sc(\mathbb{C}\times\Lambda)$.
Furthermore, we have
\[
\mathbb{H}b_\mu^{i,j}=2\pi (i+j+1) b_\mu^{i,j}
\]
\end{thm}

It follows that in $g$-dimensions, that the ground states $b_\mu$
satisfy
\[
\mathbb{H}b_\mu=2\pi g b_\mu
\]
Thus we see that there is a kernel for $\mathbb{\Delta}^0$ only
for $k-l=-g$.
\begin{thm}
\begin{enumerate} \item  The kernel of
$\mathbb{\Delta}^0$ on $\Xc^\bullet$ is zero except in dimension
$g$ where it has an orthogonal basis consisting of $\eta_\mu^0$
for $\mu\in\Lambda$ where
\begin{equation}\begin{split}
\eta_\mu^0= & b_\mu(z,\lambda)e^-_1\wedge\cdots\wedge e^-_g
\\
=   & b_\mu(z,\lambda)(d\overline{z}_1-\imath
d\overline{\zeta}_1)\wedge\cdots\wedge (d\overline{z}_g-\imath
d\overline{\zeta}_g)\end{split}\end{equation} \item The cohomology
of $(\Xc^\bullet, \mathbb{X}^0)$ is zero except in dimension $g$
where it has an orthogonal basis consisting of
$\eta_\mu^0$.\end{enumerate}
\end{thm}
Now one should note that $b_\mu$ is nothing other than
$b_0\cdot(-\mu)$. Moreover, $\mathbb{X}^0$ is linear with respect
to the right action of $\A(\Lambda;\sigma)$. Hence
\begin{cor}\label{calc-coho}
The cohomology of $(\Xc^\bullet, \mathbb{X}^0)$ is zero except in
dimension $g$ where it is a free $\A(\Lambda;\sigma)$ module of
rank one, with generator $\eta_0^0$. \end{cor}

It is now time to get the rest of $\mathbb{X}$ (not just
$\mathbb{X}^0$) into the picture. Recall
$\mathbb{X}=\mathbb{Q}\#\mathbb{P}$ and we see that
$\mathbb{X}=\mathbb{X}^0+\mathbb{X}^1$ where for $\phi\in
(\Qc^\bullet\otimes_{\A(X;B)}\Pc^\bullet)\otimes_{\A(\Lambda;\sigma)}\A^{0,\bullet}(\Lambda;\sigma)$
and using the isomorphisms described above of this with
$\Sc(V\times\Lambda;\Lambda^\bullet V_{1,0}\otimes \Lambda^\bullet
V^{0,1}\otimes\Lambda^\bullet V_{1,0})$ we have
\[
\mathbb{X}^1 \phi(z,\lambda)=(1\otimes \mathbb{P}^1)
\phi(z,\lambda)=\sum_{j=1}^g -\cf
(z_j+\lambda_j)d\overline{\tau_j}\wedge\phi(z,\lambda)
\]
Here $d\overline{\tau_j}$ is the same basis of $V_{1,0}$ as
$d\overline{\zeta_j}$ but considered in the second copy of
$V_{1,0}$. \ Now
\begin{equation}\begin{split}
 \mathbb{X}(\eta_\mu^0)&=\mathbb{X}(\eta_0^0\cdot(-\mu))\\
 & =
 \mathbb{X}(\eta_0^0)\cdot(-\mu)+\eta_0^0\dbar(-\mu)\\
 & =\sum_j -\cf
 z_j d\overline{\tau_j}\wedge\eta_0^0\cdot(-\mu)-\cf
 \eta_0^0\mu_jd\overline{\tau_j}\end{split}
\end{equation}

So none of the $\eta_\mu^0$ are closed in the complex
$(\Xc^\bullet\otimes_{\A(\Lambda;\sigma)}\A^{0,\bullet}(\Lambda;\sigma),
\mathbb{X})$.
\begin{equation}\begin{split}
\mathbb{X}(\eta_0^0)(z,\lambda)&
=\mathbb{X}^0\eta^0_0(z,\lambda)+\mathbb{X}^1\eta_0^0(z,\lambda)\\
 &=\mathbb{X}^1\eta_0^0(z,\lambda)\\
 &=\sum_j -\cf
 z_jd\overline{\tau_j}\wedge\eta_0^0(z,0)
 \end{split}
 \end{equation}
 for $\lambda=0$ and is zero for $\lambda\ne 0$.
Therefore, letting
\begin{equation}\begin{split}
\eta= & b_0(z,\lambda)(e^-_1+\imath
d\overline{\tau}_1)\wedge\cdots\wedge (e^-_g+\imath
d\overline{\tau}_g)
\\
=   & b_0(z,\lambda)(d\overline{z}_1-\imath
d\overline{\zeta}_1+\imath d\overline{\tau}_1)\wedge\cdots\wedge
(d\overline{z}_g-\imath d\overline{\zeta}_g+\imath
d\overline{\tau}_g)\end{split}\end{equation} one sees after a
somewhat lengthy but simple computation that
\[
\mathbb{X}(\eta)=0
\]
\begin{cor}\label{calc-coho2}
The cohomology of the complex
$(\Xc[g]^\bullet\otimes_{\A(\Lambda;\sigma)}\A^{0,\bullet}(\Lambda;\sigma),
\mathbb{X})$ is zero, except in dimension $0$ where it is one
(complex) dimensional.\end{cor} \proof Define a map of complexes
\[
(\A(\Lambda;\sigma)\otimes_{\A(\Lambda;\sigma)}\A^{0,\bullet}(\Lambda;\sigma),\dbar)\to
(\Xc[g]^\bullet\otimes_{\A(\Lambda;\sigma)}\A^{0,\bullet}(\Lambda;\sigma),
\mathbb{X}) \] by
\[
1\mapsto \eta
\]
Now we use the natural spectral sequence for the complex.
According to the previous corollary, \ref{calc-coho}, map above
induces an isomorphism on the $E_1$ term of this spectral
sequence. and so the cohomology identifies with
$H^*(\A^{0,\bullet}(\Lambda;\sigma))$. \qed

As a result of the two corollaries, \ref{calc-coho} and
\ref{calc-coho2} we have
\begin{prop}
$(\Xc[g]^\bullet,\mathbb{X})$ viewed as an object in
$q\Pc_{\A^{0,\bullet}(\Lambda;\sigma)}$ (i.e. just as a right
module) is quasi-isomorphic to $(\A,\dbar)$.
\end{prop}

Unfortunately, this is not quite enough to see that it induces
quasi-equivalence of categories. We need a map of quasi-perfect
twisted bimodules.

We now show that the quasi-perfect twisted bimodule
$(\Xc[g]^\bullet,\mathbb{X})$ induces a functor which is naturally
quasi-equivalent to the identity functor. It is clear that the
identity functor is implemented by the perfect twisted bimodule
$(\A(\Lambda;\sigma),\dbar)$. Using $\alpha$ from \eqref{alpha},
we define a morphism (also called $\alpha$
\begin{equation}\label{alpha1}
\alpha\in
\Hom^0_{q\Pc_{\A^{0,\bullet}(\Lambda;\sigma)}}(\Xc[g]^\bullet,\A(\Lambda;\sigma))
\end{equation}
where for $q\otimes p\in
\Xc[g]^\bullet\otimes_{\A(\Lambda;\sigma)}\A^{0,\bullet}(\Lambda;\sigma)=(\Qc^\bullet\otimes_{\A(X;B)}\Pc^\bullet)\otimes_{\A(\Lambda;\sigma)}\A^{0,\bullet}(\Lambda;\sigma)$
\begin{equation}\label{alphaqp}
\alpha(q\otimes p)= \sum_\lambda
\frac{1}{(2\sqrt{-1})^g})\left[\int_V
q(z+\lambda)p(z)\sigma^{-1}(\lambda,z)\wedge
dz_1\wedge\cdots\wedge dz_g \right][\lambda] \end{equation} Here
the $d\overline{z}$'s in the integrand combine with the $dz$'s and
are integrated. Only the top degree in $d\overline{z}$'s
contribute to the integral and thus $\alpha$ indeed maps
$\Xc[g]^\bullet\to \A(\Lambda;\sigma)$. We also note that
$d\overline{\zeta}_j$ and $d\overline{\tau}_j$ are both mapped to
$d\overline{\tau}_j$ in $\A^{0,\bullet}(\Lambda;\sigma)$. In terms
of the isomorphism of
$(\Qc^\bullet\otimes_{\A(X;B)}\Pc^\bullet)\otimes_{\A(\Lambda;\sigma)}\A^{0,\bullet}(\Lambda;\sigma)$
with $\Sc(V\times\Lambda;\Lambda^\bullet V_{1,0}\otimes
\Lambda^\bullet V^{0,1}\otimes\Lambda^\bullet V_{1,0})$ we have
for $\phi(z,\lambda)\in \Sc(V\times\Lambda;\Lambda^\bullet
V_{1,0}\otimes \Lambda^\bullet V^{0,1}\otimes\Lambda^\bullet
V_{1,0})$
\begin{equation}\label{alpha2}
\alpha(\phi)= \sum_\lambda \frac{1}{(2\sqrt{-1})^g}\left[\int_V
\phi(z+\lambda,-\lambda)\sigma^{-1}(\lambda,z)\wedge
dz_1\wedge\cdots\wedge dz_g \right][\lambda]
\end{equation}
We now show that
\begin{equation}\label{compat}
\alpha(\mathbb{X}(\phi))=\dbar(\alpha(\phi))
\end{equation}
and thus $A$ is a map of twisted bimodules. We compute

\begin{equation}\label{a1}
\begin{split}
\alpha(\mathbb{X}(\phi))&= \sum_\lambda
\frac{1}{(2\sqrt{-1})^g}\left[\int_V
\mathbb{X}(\phi)(z+\lambda,-\lambda)\sigma^{-1}(\lambda,z)\wedge
dz_1\wedge\cdots\wedge dz_g \right][\lambda]
\\
 &= \sum_{\lambda, i} \left[\frac{1}{(2\sqrt{-1})^g}\int_V \left\{\cf d\overline{\zeta}_i\wedge
(z+\lambda)_i \phi (z+\lambda,-\lambda)\right.\right.
\\
 & + d\overline{z}_i\wedge \frac{\partial
 \phi}{\partial\overline{z}_i}(z+\lambda,-\lambda) \\
  & + \cf d\overline{z}_i\wedge(\overline{B_i(-\lambda)}\phi(z+\lambda,-\lambda))\\
  &  -\cf (z+\lambda-\lambda)_i
 \left.\left. d\overline{\tau}_i\phi(z+\lambda,-\lambda)\right\}
  \sigma^{-1}(\lambda,z)\wedge
dz_1\wedge\cdots\wedge dz_g \right][\lambda]
\end{split}\end{equation}
Now we send both $d\overline{\zeta}_i$ and $d\overline{\tau}_i$ to
$d\overline{\tau}_i$
\begin{equation}\begin{split}
 & = \sum_{\lambda,i} [\frac{1}{(2\sqrt{-1})^g}\int_V \{\cf d\overline{\tau}_i\wedge
(z+\lambda)_i \phi (z+\lambda,-\lambda)
\\
&  -\cf (z+\lambda-\lambda)_i
  d\overline{\tau}_i\phi(z+\lambda,-\lambda) \\
 & + d\overline{z}_i\wedge \frac{\partial
 \phi}{\partial\overline{z}_i}(z+\lambda,-\lambda) \\
  & + \cf
  d\overline{z}_i\wedge(\overline{B_i(-\lambda)}\phi(z+\lambda,-\lambda))
  \}
  \sigma^{-1}(\lambda,z)\wedge
dz_1\wedge\cdots\wedge dz_g ][\lambda]
\end{split}
\end{equation}
which is equal to
\begin{equation}\begin{split}
 &  \sum_{\lambda,i} [\frac{1}{(2\sqrt{-1})^g}\int_V \{\cf \lambda_i d\overline{\tau}_i\wedge
 \phi (z+\lambda,-\lambda)\sigma^{-1}(\lambda,z)
\\
 & + d\overline{z}_i\wedge \frac{\partial
 \phi}{\partial\overline{z}_i}(z+\lambda,-\lambda)\sigma^{-1}(\lambda,z) \\
  & + \cf
  d\overline{z}_i\wedge(\overline{B_i(-\lambda)}\phi(z+\lambda,-\lambda))\sigma^{-1}(\lambda,z)
  \}
  \wedge
dz_1\wedge\cdots\wedge dz_g ][\lambda]
\end{split}
\end{equation}
And finally we have \begin{equation}\begin{split}
&\sum_{\lambda,i} [\frac{1}{(2\sqrt{-1})^g}\int_V \{\cf \lambda_i
d\overline{\tau}_i\wedge
 \phi (z+\lambda,-\lambda)\sigma^{-1}(\lambda,z)
\\
 & + d\overline{z}_i\wedge \frac{\partial(
 \sigma^{-1}(\lambda,z)\phi(z+\lambda,-\lambda) )}{\partial\overline{z}_i}
  \}
  \wedge
dz_1\wedge\cdots\wedge dz_g ][\lambda] \\
 & = \dbar(\alpha(\phi))+\sum_{\lambda,i} \frac{1}{(2\sqrt{-1})^g}\int_V d\overline{z}_i\wedge\frac{\partial
 (\sigma^{-1}(\lambda,z)\phi(z+\lambda,-\lambda))}{\partial
 \overline{z}_i })\wedge dz_1\wedge \cdots \wedge dz_g\\
 & =\dbar(\alpha(\phi))+\sum_{\lambda} \frac{1}{(2\sqrt{-1})^g}\int_V
 \overline{\partial}(\sigma^{-1}(\lambda,z)\phi(z+\lambda,-\lambda))\wedge dz_1\wedge \cdots \wedge dz_g\\
 & =\dbar(A(\phi))
 \end{split}
\end{equation}
This last equality holds since for any differential form $f\in
\Sc(V;\Lambda^\bullet V^{0,1})$ it follows that
\begin{equation}\begin{split}
\int_V \dbar{(f)}dz_1\wedge \cdots \wedge dz_g& =\int_V \dbar{(f
dz_1\wedge\cdots\wedge dz_g)}\\
 & =\int_V (d-\partial)(f dz_1\wedge\cdots \wedge dz_g) \\
 & =\int_V d(f dz_1\wedge\cdots\wedge dz_g)\\
 &=0
 \end{split}\end{equation}
by Stokes theorem and since $f$ is Schwartz. Thus we may apply our
criterion \eqref{lem:dgeq}  to conclude that
$(\Xc^\bullet,\mathbb{X})$ implements a functor equivalent to the
identity functor. (Since there is no curvature on this side of the
equivalence, $\Phi=0$.)

In the classical case of Mukai duality on tori, the proof that the
composition one direction gives the identity is exactly the same
calculation as the other, since both are tori. In our case, one is
a noncommutative torus and the other is a gerby torus and things
are not quite as symmetric, though as we carry the details out
below, one will see a significant overlap.

Write $\mathbb{Y}$ for the composed $\Zb$-connection
$\mathbb{P}\#\mathbb{Q}$ on the quasi-perfect twisted bimodule
$$\Yc^\bullet=\Pc^\bullet\otimes_{\A(\Lambda;\sigma)}\Qc^\bullet.$$
We calculate the zero component
$\mathbb{Y}^0=(\mathbb{P}\#\mathbb{Q})^0=\mathbb{P}\otimes
\mathbb{1}+\mathbb{1}\otimes \mathbb{Q}^0$.
\newline We have
$(\mathbb{P}\otimes \mathbb{1}+\mathbb{1}\otimes
\mathbb{Q}^0)(p\otimes q)=$
\begin{equation}\label{b1}
  \dbar_V(p)\otimes q +\cf \omega\wedge p
\otimes q -\cf D'p\otimes q +\cf p\otimes D'q
\end{equation}
Now we would like to write down more explicitly quasi-perfect
twisted bimodule $(\Yc^\bullet,\mathbb{Y})$.
\begin{equation}\begin{split}
\Yc^\bullet & \equiv \Pc^\bullet\otimes _{\A(\Lambda;\sigma)}\Qc^\bullet\\
 & \cong \Sc(V;\Lambda^\bullet V^{0,1})\otimes
 _{\A(\Lambda;\sigma)}\Sc(V;\Lambda^\bullet V_{1,0})
 \end{split}\end{equation}
This last expression is the quotient of $\Sc(V;\Lambda^\bullet
V^{0,1})\otimes
 _\mathbb{C}\Sc(V;\Lambda^\bullet V_{1,0})$ by the closure of the relation
 $p\lambda\otimes q=p\otimes \lambda q$
 or what is the same thing, $p\lambda\otimes\lambda^{-1}q=p\otimes q$. This is the coinvariants
 by the right action of  $\Lambda$ on
 \[
 \Sc(V\times V;\Lambda^\bullet V^{0,1}\otimes
 \Lambda^\bullet V_{1,0})
 \]
given by \[ (\phi\cdot
\lambda)(z,w)=\phi(z+\lambda,w+\lambda)\sigma(\lambda,z-w)
 \]
where $\phi\in \Sc(V\times V;\Lambda^\bullet V^{0,1}\otimes
 \Lambda^\bullet V_{1,0})$, $\lambda\in\Lambda$. Let us emphasize
 that here by the coinvariants we mean
 \[
\Sc(V\times V;\Lambda^\bullet V^{0,1}\otimes
 \Lambda^\bullet V_{1,0}) / (\mbox{closure of span }(\phi-\phi\cdot\lambda))
 \]
 \begin{prop}
\begin{enumerate}
\item $\Yc^\bullet\cong (\Sc(V\times V;\Lambda^\bullet
V^{0,1}\otimes
 \Lambda^\bullet V_{1,0}))_\Lambda $
 \[
\cong \{\phi\in C^\infty(V\times V; \Lambda^\bullet
V^{0,1}\otimes\Lambda^\bullet V_{1,0})^\Lambda | \,\, \phi \mbox{
satisfies the Schwartz estimates in } z-w \}.
 \]
That is, $\phi$ is invariant and satisfies
\[
(z-w)^\alpha\frac{\partial^{\beta,\gamma}\phi}{\partial z^\beta
w^\gamma}\in L^\infty(V\times V;\Lambda^\bullet
V^{0,1}\otimes\Lambda^\bullet V_{1,0})
\]
for all multi-indices $\alpha$, $\beta$ and $\gamma$.

\item Under this isomorphism, $\mathbb{Y}^0$ is
\[
\mathbb{Y}^0\phi(z,w)= \overline{\partial}_z\phi(z,w)
 +\cf \omega(z)\wedge \phi -\cf D'(z)\phi(z,w) +\cf D'(w)\phi(z,w)
\]
and using the same conventions as before (in particular,
\eqref{def:B_j}) this can be expressed as
\[
\mathbb{Y}^0\phi(z,w)=\sum_j
d\overline{z_j}\wedge\frac{\partial\phi}{\partial \overline
{z_j}}(z,w)
 +\cf d\overline{z_j}\wedge\overline{B_j(z)}\phi(z,w) +\cf d\overline{\zeta_j}\wedge(w_j-z_j)\phi(z,w)
\]
\end{enumerate}
 \end{prop}
 \proof

Let us call the space of invariants described in the proposition
$W$.
 The map implementing the isomorphism is $\tau:(\Sc(V\times V;\Lambda^\bullet
V^{0,1}\otimes
 \Lambda^\bullet V_{1,0}))_\Lambda\to W$ is
 \[
\tau(\phi)(z,w) = \sum_\lambda
\phi(z+\lambda,w+\lambda)\sigma(\lambda,z-w)
 \]
Clearly this map is well defined on the coinvariants and injective
and one checks that the image is in $W$.  One defines a section
$\rho:W\to (\Sc(V\times V;\Lambda^\bullet V^{0,1}\otimes
 \Lambda^\bullet V_{1,0}))$ by
\[
\rho(\psi)(z,w)= h(z)\psi(z,w)
\]
where $h$ is a function as in the proof of \ref{proj-lemma}. The
image of this in the coinvariants is a section.

The computation of $\mathbb{Y}$ under this isomorphism is clear.
\qed.

\begin{rem} It seems to be a general phenomenon that the space of
 coinvariants of a space of functions under a proper group action
 can be expressed as a space of invariants. One can also write
 $W$ as the sections of an infinite dimensional vector bundle over
 $X$. From now on, when we talk about $\Yc^\bullet$, we will
 implicitly use the isomorphism with $W$.
 \end{rem}

We will be applying the criterion\eqref{lem:dgeq} to conclude that
$\Yc_*$ is naturally DG-quasi-equivalent to the identity functor.
Unlike the situation for $\mathbb{X}$, our connection $\mathbb{Y}$
has curvature. In particular $(\mathbb{Y}^0)^2=\cf B^{0,2}$ We
define an endomorphism \[ \Phi:\Yc^\bullet\to\Yc^\bullet
\] of degree one by
\[
\Phi(\phi)(z,w)=\cf\sum_j
\overline{B_j(w-z)}d\overline{z}_j\wedge\phi(z,w)
\]
Then $\Phi\circ\Phi=0$ and we have $[\mathbb{Y}^0,\Phi]=-\cf
B^{0,2}$. Now we need to calculate \[H^*(\Yc^\bullet,
\mathbb{Y}^0+\Phi).\]

We have $(\mathbb{Y}^0+\Phi)\phi(z,w)$
\[
=\sum_j d\overline{z_j}\wedge\frac{\partial\phi}{\partial
\overline {z_j}}(z,w)
 +\cf d\overline{z_j}\wedge\overline{B_j(w)}\phi(z,w) +\cf d\overline{\zeta_j}\wedge(w_j-z_j)\phi(z,w)
\]

We calculate the Laplacian of $\mathbb{Y}^0+\Phi$. The adjoint
\[ (\mathbb{Y}^0+\Phi)^*(\phi)(z,w)=\sum_j -\frac{\partial}{\partial
z_j}\iota_{\frac{\partial}{\partial \overline{z_j}}}\phi(z,w) -\cf
B_j(w)\iota_{\frac{\partial}{\partial
\overline{z_j}}}\phi(z,w)-\cf
\overline{(w_j-z_j)}\iota_{\frac{\partial}{\partial
\overline{\zeta_j}}}\phi(z,w)
\]
Write $Y_j(w)$ for the first order differential operator
\[
Y_j(w)(\phi)(z,w)=\overline{B_j(w)}\frac{\partial}{\partial
z_j}\phi+ B_j(w)\frac{\partial}{\partial \overline{z_j}}\phi(z,w)
\]
Then the Laplacian
$\bbb^0(\phi)(z,w)=\left((\mathbb{Y}^0+\Phi)^*(\mathbb{Y}^0+\Phi)+(\mathbb{Y}^0+\Phi)(\mathbb{Y}^0+\Phi)^*\right)(\phi)(z,w)$
\begin{equation}\label{e1}
 =
\sum_{j=1}^g\left(-\frac{\partial^2}{\partial z_j\partial
\overline{z_j}}- \cf Y_j(w)+4\pi^2(|z_j-w_j|^2+|B_j(w)|^2)+\cf
(d\overline{z}_j\circ
\iota_{\frac{\partial}{\partial\overline{\zeta}_j}}+\iota_{\frac{\partial}
{\partial\overline{z}_j}}\circ d\overline{\zeta}_j)\right)
\phi(z,w)\end{equation}\[
=\sum_{j=1}^g\left(-\frac{\partial^2}{\partial z_j\partial
\overline{z_j}}- \cf Y_j(w)+4\pi^2(|z_j-w_j|^2+|B_j(w)|^2)+\cf
L_j)\right)\phi(z,w)
\]
where we recall that $L_j=d\overline{z}_j\circ
\iota_{\frac{\partial}{\partial\overline{\zeta}_j}}+\iota_{\frac{\partial}
{\partial\overline{z}_j}}\circ d\overline{\zeta}_j$ and that we
have for $I=(i_1< i_2< \cdots < i_k)$ and $J=(j_1 < j_2< \cdots <
j_l)$
\[
L(e^+_I\wedge e^-_J)=(k-l)\imath(e^+_I\wedge e^-_J)
\]
where $L=\sum_j L_j$ and
\[
e_j^\pm=d\overline{z}_j\pm \imath d\overline{\zeta}_j
\]
Hence we have the eigenvector decomposition for $L$
\[
\Lambda^\bullet V_{1,0}\otimes \Lambda^\bullet
V^{0,1}=\oplus_{I,J}\mbox{ span }e^+_I\wedge e^-_J
\]
Thus we have that on the $e^+_I\wedge e^-_J$ component,
$\bbb^0\phi (z,w)$
\begin{equation}\label{b1}
= \left\{\sum_j\left(-\frac{\partial^2}{\partial z_j\partial
\overline{z_j}}- \cf
Y_j(w)+4\pi^2(|z_j-w_j|^2+|B_j(w)|^2)\right)-2\pi(k-l)\right\}
\phi(z,w).
\end{equation}

Now we solve this equation. First, note that in the definition of
the $2$-cocycle $\sigma$, we may insert elements from $V$ into
both arguments. Thus $\sigma(z,w)$ is a well defined function on
$V\times V$.  One has the following formulae:
\begin{equation}\label{eq:sigma-dif}\begin{split}\frac{\partial}{d\overline{z}_j}\sigma(z,w)&=-\cf
\overline{B_j(w)}\sigma(z,w)\\
\frac{\partial}{dz_j}\sigma(z,w)&=-\cf
B_j(w)\sigma(z,w)\end{split}
\end{equation}
Define
\begin{equation}\label{def:vacstate}
a(z,w) =\sigma(z,w)\exp(-2\pi(|z-w|^2)
\end{equation}
Then
\begin{equation}\begin{split}
a(z+\lambda,w+\lambda)\sigma(\lambda,z-w)&=
\sigma(z+\lambda,w+\lambda)\exp(-2\pi|(z+\lambda)-(w+\lambda)|^2)\sigma(\lambda,z-w)\\
 &
 =\sigma(\lambda,w)\sigma(z,\lambda)\sigma(z,w)\exp(-2\pi|z-w|^2)\sigma(\lambda,z-w)\\
  &= a(z,w)
  \end{split}\end{equation}
  so $a\in \Yc^0$.

Let's write $\mathbb{L}$ for the operator
\[
\mathbb{L}\phi(z,w)=\sum_j\left(-\frac{\partial^2}{\partial
z_j\partial \overline{z_j}}- \cf
Y_j(w)+4\pi^2(|z_j-w_j|^2+|B_j(w)|^2)\right)\phi(z,w) \]  defined
on $\Yc^0$.

We now show that $a$ is an eigenvector for $\mathbb{L}$. We have
from \eqref{eq:sigma-dif} the following formulae
\begin{enumerate}
\item $\frac{\partial a}{\partial\overline{z}_j}(z,w)  =
-2\pi((z_j-w_j)+\imath \,\,\overline{B_j(w)})a(z,w)$ \item
$\frac{\partial a}{\partial z_j}(z,w)  =
-2\pi(\overline{(z_j-w_j)}+\imath B_j(w))a(z,w)$ \item
$\frac{\partial^2 a}{\partial z_j\partial\overline{z_j}}(z,w) =
-2\pi a(z,w)\\
+4\pi^2\left(|z_j-w_j|^2+\imath(z_j-w_j)B_j(w)+\imath\,\,\overline{(z_j-w_j)B_j(w)}-
|B_j(w)|^2\right)a(z,w)$ \item $ Y_j(z)a(z,w)  = -2\pi B_j(w)
((z_j-w_j)+\imath\,\,\overline{B_j(w)})a(z,w)\\
  + -2\pi\overline{B_j(w)}
(\overline{(z_j-w_j)}+\imath\,\,B_j(w))a(z,w)$
\end{enumerate}
So
\begin{equation}\begin{split} \mathbb{L}a(z,w) & =\sum_j 2\pi
a(z,w)-\sum_j 4\pi^2|z_j-w_j|^2\\
 &-\sum_j 4\pi^2\imath\left(B_j(w)(z_j-w_j)+\,\,\overline{B_j(w)}\overline{(z_j-w_j)}\right)a(z,w)\\
 &+\sum_j 4\pi^2\left(|B_j(w)|^2\right)a(z,w)\\
  & +\sum_j
  4\pi^2\imath(B_j(w)((z_j-w_j)+\imath \,\,\,\,\overline{B_j(w)})\\
  & +4\pi^2\imath\,\,\,\,\overline{B_j(w)}(\overline{(z_j-w_j)}+\imath B_j(w))a(z,w)\\
   & +\sum_j 4\pi^2(|z_j-w_j|^2+|B_j(x)|^2)a\\
   & = 2\pi g a(z,w)\\
\end{split}\end{equation}

We calculate the eigenvalues and eigenvectors of $\mathbb{L}$. As
before, we do this in the one dimensional case.  Define operators
\begin{equation}\label{anncreops2}
\begin{split} \mathbb{A}\phi(z,w) & =\left(\frac{\partial}{\partial
\overline{z}}+2\pi((z-w)+\imath B(w))\right)\phi\\
 \mathbb{A}^*\phi(z,w) & =\left(-\frac{\partial}{\partial
z}+2\pi(\overline{(z-w)}-\imath \overline{B(w)})\right)\phi\\
 \mathbb{B}\phi (z,w)& =\left(\frac{\partial}{\partial
z}+2\pi(\overline{(z-w)}+\imath B(w))\right)\phi \\
 \mathbb{B}^*\phi(z,w) & =\left(-\frac{\partial}{\partial
\overline{z}}+2\pi((z-w)-\imath \overline{B(w)})\right)\phi
\end{split}\end{equation}
These operators satisfy the following commutation relations.
\begin{equation}\label{commeq}
\begin{split}
[\mathbb{L},\mathbb{A}]=-2\pi\mathbb{A}, & \,\,\,\,
[\mathbb{L},\mathbb{A}^*]=2\pi\mathbb{A}^*\\
[\mathbb{L},\mathbb{B}]=-2\pi\mathbb{B}, &\,\,\,\,
[\mathbb{L},\mathbb{B}^*]=2\pi\mathbb{B}^*\\
[\mathbb{A},\mathbb{A}^*]=4\pi, &\,\,\,\,
[\mathbb{B},\mathbb{B}^*]=4\pi \\
[\mathbb{A},\mathbb{B}]=0, & \,\,\,\,[\mathbb{A},\mathbb{B}^*]=0\\
[\mathbb{A}^*,\mathbb{B}]=0,  &\,\,\,\,
[\mathbb{A}^*,\mathbb{B}^*]=0
\end{split}\end{equation}
Set $a^{0,0}(z,w)=a(z,w)$ from above. Define recursively,
\begin{equation}\label{Hermite}
\begin{split}
a^{i+1,j}(z,w)=\mathbb{A}^*a^{i,j}(z,w), & \,\,\,\,
a^{i,j+1}(z,w)=\mathbb{B}^*a^{i,j}(z,w)
\end{split}
\end{equation}
This is well defined since $\mathbb{A}^*$ and $\mathbb{B}^*$
commutate. Moreover, $\mathbb{A}a^{0,0}=\mathbb{B}a^{0,0}=0$. Then
as before we have
\begin{thm}
For each $w\in V$, the functions $a^{i,j}(\cdot,w)\in
\Sc(\mathbb{C})$ form an orthogonal complete basis of the closure
$L^2(\mathbb{C})$ of $\Sc(\mathbb{C})$. Furthermore, we have
\[
\mathbb{L}a^{i,j}=2\pi (i+j+1) a^{i,j}
\]
\end{thm}

It follows that in $g$-dimensions, that the ground states $a$
satisfy
\[
\mathbb{L}a=2\pi g a
\]
Thus we see that there is a kernel for $\mathbb{\bbb}^0$ only for
$k-l=g$.
\begin{thm}\label{calc-coho-Y}
The cohomology of $(\Yc^\bullet, \mathbb{Y}^0+\Phi)$ is zero
except in dimension $g$ where it is a free $\A(\Lambda;\sigma)$
module of rank one, with generator $a$. \end{thm}

Finally, define a map of $\A^{0,\bullet}(X;B)$-bimodules
\[
\beta:\Yc^\bullet\otimes_{\A(X;B)}\A^{0,\bullet}(X;B)\to
\A^{0,\bullet}(X;B)
\]
by
\[
\beta(\phi)(z)=\iota_\Xi(\phi)(z,z)
\]
As before, this needs a little explaining. We have that
\[\Yc^\bullet\otimes_{\A(X;B)}\A^{0,\bullet}(X;B)\cong
\Yc^\bullet\otimes \Lambda^\bullet V^{0,1}.\] In this last factor
of $\Lambda^\bullet V^{0,1}$ we denote the basis by
$d\overline{w}_j$. Also, $\Xi$ is the alternating multivector
$\frac{\partial}{\partial\overline{\zeta}_1}\wedge\cdots\wedge
\frac{\partial}{\partial\overline{\zeta}_g}$. Then the map $\beta$
does the following. It picks off any factor containing
$d\overline{\zeta}_1\wedge\cdots\wedge d\overline{d\zeta}_g$, it
sends both $d\overline{z}_j$ and $d\overline{w}_j$ to
$d\overline{z}_j$ in $\A^{0,\bullet}(X;B)$ and it restricts this
to the diagonal.
\begin{prop}
The map $\beta$ is a map for $\A^{0,\bullet}(X;B)$-bimodules and
it commutes with the $\mathbb{Z}$-connections.
\end{prop}

This completes the proof of theorem \ref{theorem-main}.

\end{document}